

\input amstex
\magnification=1100
\voffset=-.25truein
\hoffset=3truepc
\magnification=1100

\def\ie{i\rom.e\rom.}
\def\eg{e\rom.g\rom.}

\def\kl{}
\def\Kl{^{\mathstrut}} 
\def\itemitemitem{\par\indent\indent\hangindent3\parindent\textindent}
\loadeusm

\input epsf.tex 
\documentstyle{amsppt}
\topmatter

\title
Diagrams with selection
and method for constructing
boundedly generated
and boundedly simple groups
\endtitle

\author
Alexey Muranov
\endauthor

\affil
Vanderbilt University
\endaffil

\address
Vanderbilt University, Dept of Mathematics\newline
\indent 1326 Stevenson Center\newline
\indent Nashville TN 37240--0001
\endaddress

\email
alexey.muranov\@vanderbilt.edu
\endemail

\date
April 26, 2004
\enddate


\thanks
This work was supported in part by the NSF grant DMS 0072307
of Alexander Ol'shanskii and Mark Sapir.
\endthanks



\abstract\nofrills
The existence of an infinite simple boundedly generated $2$-generated group
and the existence of a boundedly simple $2$-generated group
containing a free non-cyclic subgroup are proved.
\endabstract

\toc
\widestnumber\head{14.}
\head 1. Introduction\endhead
\head 2. Graphs and maps\endhead
\head 3. Maps with selection\endhead
\head 4. Condition $\Cal A$\endhead
\head 5. Estimating lemmas\endhead
\head 6. Statement of The Main Theorem\endhead
\head 7. Inductive Lemma 1\endhead
\head 8. Inductive Lemma 2\endhead
\head 9. Proof of The Main Theorem\endhead
\head 10. Lemma about exposed face\endhead
\head 11. Group presentations and diagrams\endhead
\head 12. Proof of Theorem 1\endhead
\head 13. Proof of Theorem 2\endhead
\endtoc

\rightheadtext{Diagrams with selection}
\endtopmatter
\document

\head
1. Introduction
\endhead

\definition{Definition}
A group $G$ is called {\it $m$-boundedly generated\/} if
it has (not necessarily normal) cyclic subgroups
$C_1\kl$, \dots, $C_m\kl$ such that $G=C_1\kl\cdots C_m\kl$.
A group $G$ is called {\it boundedly generated\/} if
it is $m$-boundedly generated for some natural~$m$.
\enddefinition

Apparently, the definition of bounded generation was motivated by the work
\cite{CK} of Carter and Keller, in which they proved that any matrix in
$SL_n(\Cal O)$, where $n\ge3$ and $\Cal O$ is the ring of integers of
a finite extension of the field of rational numbers,
is the product of a bounded number of elementary
matrices (bounded depending only on $n$ and~$\Cal O$).
Since then, bounded generation has been studied in connection with
the {\it congruence property\/ $($CSP\/$)$} (see \cite{PR}),
{\it Kazhdan's property\/ $($T\/$)$} (see~\cite{Sh}).
The following closely related property of profinite groups is used
in the book \cite{LuSe} by Lubotzky and Segal, as well as in \cite{PR}:
a profinite group is {\it boundedly generated\/} as a profinite group if
it is the product of finitely many procyclic subgroups.

\definition{Definition}
A group $G$ is called {\it $m$-boundedly simple\/} if
for any two nontrivial elements $g,h\in G$,
the element $h$ may be presented as the product of $m$ or fewer
conjugates of~$g^{\pm1}$.
A group $G$ is called {\it boundedly simple\/} if
it is $m$-boundedly simple for some natural~$m$.
\enddefinition

\remark{Remark\/ \rm1.1}
For each natural $m$, the class of $m$-boundedly simple groups is
definable by a formula of the restricted predicate calculus.
The class of all simple groups is not.
\endremark

It is obvious that a boundedly simple group is simple.
However, not every simple group is boundedly simple.
For example, an infinite alternating group is simple
but not boundedly simple.

Every non-abelian simple group coincides with its derived subgroup.
In a non-abelian $m$-boundedly simple group, the derived subgroup in addition has
finite commutator width:
each element of the derived subgroup may be
presented as the product of $m$ or fewer commutators
(take an arbitrary nontrivial commutator, every element is the product
of $m$ or fewer conjugates of this commutator).

It is not hard to prove the following

\proclaim{Proposition 1.1}
A group is boundedly simple if and only if
each of its ultrapowers is simple$.$
If a group is\/ $m$-boundedly simple$,$ then all its ultrapowers are\/
$m$-boundedly simple$.$
\endproclaim

It is not obvious that infinite boundedly simple finitely generated
groups or infinite simple boundedly generated groups exist.
The following theorem may be found in \cite{O}:

\proclaim{Theorem (S.~Ivanov)}
For every big enough prime\/ $p,$
there exists a\/ $2$-generated infinite group of exponent\/ $p,$
in which there are exactly\/ $p$ distinct conjugacy classes$,$
and therefore$,$ every subgroup of order\/ $p$
has elements from all of these classes$.$
\endproclaim

It immediately follows that the group whose existence is stated in this theorem
is $(p-1)$-boundedly simple.
The original proof of Ivanov's Theorem works only for big $p$, say, $p\ge10^{78}$.
It heavily uses techniques of graded diagrams.
The following two theorems are proved in this paper in a different manner.

\proclaim{Theorem 1}
Let\/ $F$ be a finite-rank free group with a fixed basis\rom.
Let\/ $S$ be an infinite subset of\/ $F$ whose all elements 
are cyclically reduced and are not proper powers\rom.
If the symmetrization of\/ $S$ 
satisfies the small cancellation condition\/ $C'(\lambda)$
for some\/ $\lambda<\frac{1}{13}$\rom,
then there exist an infinite simple\/ $2$-generated group\/ $P$ and
a homomorphism\/ $\phi\:F\to P$ such that\/ $\phi$ maps\/ $S$ onto\/~$P$\rom.
\endproclaim

\proclaim{Corollary}
There exists an infinite simple\/ $2$-generated group\/ $G$ and\/
$27$ elements\/ $x_1\kl,\,\ldots,x_{27}\kl\in G$ such that for every element\/ $g$ of\/ $G$
there exists a natural number\/ $n$ such that\/ $g=x_1^n\cdots x_{27}^n$\rom.
Such a group is\rom, in particular\rom, $27$-boundedly generated\rom.
\endproclaim

\proclaim{Theorem 2}
There exists a $14$-boundedly simple\/ $2$-generated group
with solvable word problem$,$ containing a free non-cyclic subgroup$.$
\endproclaim

The purpose of this paper is to develop a method for constructing groups
with these kinds of properties.
Roughly speaking, the method consists in constructing group presentations
by imposing all necessary relations, while
making sure that ``big parts'' of these relations satisfy
certain small-cancellation-type conditions.
No classical small cancellation conditions are required
to hold for the set of imposed relations itself.
The main tool for proving that the obtained groups
are nontrivial, or even contain free non-cyclic subgroups,
is properties of van Kampen {\it diagrams with selection}
that satisfy a certain condition $\Cal A$ with several parameters.
Diagrams with selection, also called here {\it S-diagrams}, are
defined and used similarly to maps with partitioned boundaries of cells
in \cite{O} and~\cite{Sa}.
They are diagrams in which some subpaths of the contours of the faces are selected.
When choosing a selection on a van Kampen diagram over a constructed
group presentation,
the subpaths corresponding to the aforementioned ``big parts''
of the relations are selected.
Among other things, the condition $\Cal A$ requires that,
with few exceptions, each arc between two faces that lies
on selected subpaths of the contours of both faces has to be ``short''
relative to the lengths of the contours of both (!) faces.
(This makes the condition $\Cal A$ essentially
different from the condition $B$ in \cite{O}.)
A condition similar to the condition $\Cal A$ was used by A.~Sanin
in his paper \cite{Sa} in~1997.
The main theorem proved in this paper implies in particular that
if a reduced disc S-diagram satisfies the condition
$\Cal A$ with values of the parameters satisfying certain inequalities,
and the diagram either has ``many'' faces, or has faces with ``long''
contours, then the contour of the diagram itself is ``long.''
This main theorem makes the proofs of Theorems~1 and~2 relatively simple.
\par

\head
2. Graphs and maps
\endhead

Some terms pertaining to graphs and $2$-complexes
shall be clarified in this section.

If $X$ is a set, then $\|X\|$ shall denote the cardinality of~$X$.

The term {\it graph\/} will refer to an object that in graph theory is called
an undirected pseudograph (multiple edges and loops are admissible).
If $\Gamma$ is a graph, then $\Gamma(0)$ denotes its set of
vertices, and $\Gamma(1)$ denotes its set of edges;
$\Gamma(1)$ may be empty, but $\Gamma(0)$ may not.
The graph itself is an ordered pair:
$\Gamma={(}\Gamma(0),\Gamma(1){)}$.
In the context of this paper, {\it graph\/} is a synonym of {\it $1$-complex}.
A great part of this paper deals with {\it $2$-complexes}.
If $\Phi$ is a $2$-complex, then $\Phi(0)$, $\Phi(1)$, and $\Phi(2)$
denote its set of vertices, its set of edges,
and its set of faces (cells) respectively.
The $2$-complex $\Phi$ itself is the ordered triple
${(}\Phi(0),\Phi(1),\Phi(2){)}$.
Below is a formal definition of combinatorial
$0$-, $1$-, and $2$-complexes and their morphisms.

A (combinatorial) {\it$0$-complex\/} $A$ is a
$1$-tuple ${(} A(0){)}$
where $A(0)$ is an arbitrary nonempty set.
Elements of $A(0)$ are called {\it vertices\/} of the complex~$A$.

A $0$-complex with exactly $2$ vertices will be called
a {\it combinatorial\/ $0$-sphere}.

A {\it morphism\/} $\phi$ of a $0$-complex $A$ to a $0$-complex $B$ is a
$1$-tuple ${(}\phi(0){)}$ where $\phi(0)$
is an arbitrary function $A(0)\to B(0)$.
If $A$, $B$, and $C$ are $0$-complexes, and $\phi\:A\to B$ and $\psi\:B\to C$ are morphisms, then
the product $\psi\phi\:A\to C$ is defined naturally:
$(\psi\phi)(0)=\psi(0)\circ\phi(0)$.
A morphism $\phi\:A\to B$ is called an {\it isomorphism\/}
of $A$ with $B$ if
there exists a morphism $\psi\:B\to A$
such that $\psi\phi$ is the {\it identity morphism\/}
of the complex $A$ and $\phi\psi$ is the identity morphism
of the complex~$B$.

A (combinatorial) {\it$1$-complex\/} $A$ is a $2$-tuple (\ie, ordered pair)
${(} A(0),A(1){)}$ such that:
\medskip

\itemitem{$(1)$}
$A(0)$ is an arbitrary nonempty set;

\itemitem{$(2)$}
$A(1)$ is a set of ordered pairs of the form ${(} E,\alpha{)}$
where $E$ is a combinatorial $0$-sphere and $\alpha$
is a morphism of $E$ to the $0$-complex ${(} A(0){)}$.
\medskip

\noindent
Elements of $A(0)$ are called {\it vertices\/} of the complex $A$,
elements of $A(1)$ are called {\it edges}.
The $0$-complex ${(} A(0){)}$ is called the {\it$0$-skeleton\/} of~$A$.
All $1$-complexes will also be called {\it graphs}.

A $1$-complex representing a circle
(\eg, a $1$-complex consisting of one vertex and one edge)
will be called a {\it combinatorial\/ $1$-sphere},
or a {\it combinatorial circle}.
(The definition of a $1$-sphere could be made precise,
but then it would become unreasonably long.)

If $A$ and $B$ are $1$-complexes, then a {\it morphism\/} $\phi\:A\to B$
is a $2$-tuple ${(}\phi(0),\phi(1){)}$ such that:
\medskip

\itemitem{$(1)$}
$\phi^0={(}\phi(0){)}$ is a morphism of the $0$-skeleton
$A^0$ of $A$ to the $0$-skeleton $B^0$ of~$B$;

\itemitem{$(2)$}
$\phi(1)$ is a function on $A(1)$ such that
the image of every $e={(} E,\alpha{)}\in A(1)$ under $\phi(1)$
is an ordered pair ${(} e',\xi{)}$ where
$e'={(} E',\alpha'{)}\in B(1)$, $\xi$
is an isomorphism of $E$ with $E'$,
and $\phi^0\alpha=\alpha'\xi$.
\medskip

\noindent
Multiplication of morphisms
of $1$-complexes is defined naturally.
For example, if $A$, $B$, $C$ are $1$-complexes,
$\phi$ and $\psi$ are morphisms,
$\phi\:A\to B$, $\psi\:B\to C$,
$e={(} E,\alpha{)}$ is an edge of $A$,
$e'={(} E',\alpha'{)}$ is an edge of $B$,
$e''={(} E'',\alpha''{)}$ is an edge of $C$,
$\phi(1)(e)={(} e',\xi{)}$,
$\psi(1)(e')={(} e'',\zeta{)}$,
then $(\psi\phi)(1)(e)={(} e'',\zeta\xi{)}$
(note that $\zeta\xi$ is an isomorphism of the $0$-complex $E$
with the $0$-complex~$E''$).
{\it Isomorphisms\/} of $1$-complexes are defined in the usual way.

A (combinatorial) {\it $2$-complex\/} $A$ is a $3$-tuple
${(} A(0),A(1),A(2){)}$ such that:
\medskip

\itemitem{$(1)$}
${(} A(0),A(1){)}$ is a $1$-complex, called the {\it$1$-skeleton\/} of~$A$;

\itemitem{$(2)$}
$A(2)$ is a set of ordered pairs of the form ${(} F,\alpha{)}$
where $F$ is a combinatorial $1$-sphere and $\alpha$
is a morphism of $F$ to the $1$-skeleton of~$A$.
\medskip

\noindent
Elements of $A(0)$ are called {\it vertices\/} of the complex $A$,
elements of $A(1)$ are called {\it edges},
elements of $A(2)$ are called {\it faces}.

If $A$ and $B$ are $2$-complexes, then a {\it morphism\/} $\phi\:A\to B$
is defined as a $3$-tuple ${(}\phi(0),\phi(1),\phi(2){)}$ such that:
\medskip

\itemitem{$(1)$}
$\phi^1={(}\phi(0),\phi(1){)}$ is a morphism of the $1$-skeleton
$A^1$ of $A$ to the $1$-skeleton $B^1$ of~$B$;

\itemitem{$(2)$}
$\phi(2)$ is a function on $A(2)$ such that
the image of every $f={(} F,\beta{)}\in A(2)$ under $\phi(2)$
is an ordered pair ${(} f',\xi{)}$ where
$f'={(} F',\beta'{)}\in B(2)$, $\xi$
is an isomorphism of $F$ with $F'$,
and $\phi^1\beta=\beta'\xi$.
\medskip

\noindent
Products of morphisms of $2$-complexes are defined analogously
to the case of $1$-complexes.
The notion of {\it isomorphism\/} for $2$-complexes is the natural one.

A combinatorial $n$-complex $A$ is called {\it finite\/} if
all of the sets $A(0)$, \dots, $A(n)$ are finite.
An $n$-complex $A$ is a {\it subcomplex\/} of an $m$-complex $B$,
$0\le n\le m$, if
$A(0)\subset B(0)$, \dots, $A(n)\subset B(n)$.
Every morphism $\phi$ of a combinatorial $n$-complex $A$
to a combinatorial $n$-complex $B$, $n\le 2$,
naturally defines functions
$\bar\phi^i\:A(i)\to B(i)$, $0\le i\le n$.
This notation for these functions associated with a given
morphism $\phi$ shall be used in this section for brevity.

If $e={(} E,\alpha{)}$ is an edge of a graph $\Gamma$,
then the vertices of $\Gamma$ that are the images
of the vertices of $E$ under $\bar\alpha^0$
are called the {\it end-vertices\/} of~$e$.
An edge is {\it incident\/} to its end-vertices.
A {\it loop\/} is an edge that has only one end-vertex.
Two vertices are called {\it adjacent\/} if
they form the set of end-vertices of some edge.
The end-vertex of a loop is adjacent to itself.
If $v$ is a vertex of a graph $\Gamma$, then
the number of all edges of $\Gamma$ incident to $v$ plus the number
of all loops incident to $v$ is called the {\it degree\/}
of the vertex $v$ and is denoted by $d(v)$ or~$d_\Gamma\kl(v)$.

If $f={(} F,\beta{)}$ is a face of a $2$-complex $\Phi$,
then $f$ is said to be {\it incident\/}
to the images of the vertices and edges of $F$
under $\bar\beta^0$ and $\bar\beta^1$ respectively.
If the combinatorial circle $F$ has $n$ edges ($\|F(1)\|=n\in\Bbb N$),
then the number $n$ is called the {\it degree\/} of $f$
and is denoted by $d(f)$ or~$d_\Phi\kl(f)$.

It is somewhat complicated to talk about combinatorial
complexes in purely combinatorial terms.
The geometrical intuition may help if together with
every combinatorial complex consider some corresponding topological space.

Given a $2$-complex $\Phi$, put
a $1$-point topological space $D_v^0$ into correspondence
to every vertex $v$ of $\Phi$,
a topological closed segment $D_e^1$ to every edge $e$ of $\Phi$,
and a topological closed disk $D_f^2$ to every face $f$ of~$\Phi$.
Assume that $D_x^m$ and $D_y^n$ are disjoint unless
$m=n$ and $x=y$.
Consider the topological sum
$$
\sum_{v\in\Phi(0)}D_v^0+\sum_{e\in\Phi(1)}D_e^1+\sum_{f\in\Phi(2)}D_f^2.
$$
Take the quotient of it over an equivalence relation in accordance with
the structure of~$\Phi$.
For example, if $\{v_1\kl,v_2\kl\}$ is the set of end-vertices
of an edge $e$ of $\Phi$, then one of the end-vertices
of $D_e^1$ should be identified with the only point of $D_{v_1\kl}^0$,
and the other end-vertex should be identified with the only point of~$D_{v_2\kl}^0$.
Taking the quotient may be done in two steps:
first, attach the end-vertices of the segments
$\{D_e^1\}_{e\in\Phi(1)}\kl$
to the points of the discrete topological space
$\sum_{v\in\Phi(0)}D_v^0$;
then, attach the boundaries of the discs
$\{D_f^2\}_{f\in\Phi(2)}\kl$ to the obtained ``skeleton.''
The constructed topological quotient space is defined uniquely
up to homeomorphism.
This space or any homeomorphic one is called
a {\it topological space of the complex\/}~$\Phi$.
It may also be said that $\Phi$
{\it represents\/} this topological space.
Note that the restriction of the quotient
function to the interior (in the geometric sense)
of every $D_e^1$, $e\in\Phi(1)$, and every $D_f^2$, $f\in\Phi(2)$,
is a homeomorphism onto the image in the quotient space.
Every morphism of combinatorial complexes defines
some set of continuous functions from
a given topological space of the first complex
to a given topological space of the second.

From now on, combinatorial and topological languages
shall be used together.
Moreover, some facts intuitively clear from topological point of view
shall be used without proofs.

It is a simple but tiresome task to formulate a combinatorial
criterion in terms of the local structure
of a $2$-complex ({\it stars\/} at its vertices)
that would determine if a topological space of this complex
is a $2$-dimensional surface
(\ie, a $2$-manifold or $2$-manifold with a boundary).
For an example of constructing surfaces combinatorially
using simplicial complexes see~\cite{A}.
Keeping this in mind, the following definition
may be viewed as combinatorial.
A {\it combinatorial surface\/} is a $2$-complex representing
some surface.
A {\it combinatorial sphere\/} and a {\it combinatorial disc\/} are
$2$-complexes representing a sphere and a disc respectively.
They play an important role in this paper.
Every combinatorial sphere and every combinatorial disc
are finite (because of the compactness of spheres and discs).
Combinatorial spheres and combinatorial discs may also be
defined combinatorially.

Now, a definition of a planar graph may be given easily.
A finite graph $\Gamma$ is {\it planar\/} if
it is a subgraph of
the $1$-skeleton of some combinatorial sphere.
Actually, every {\it connected\/} finite planar graph having at least
one edge is the $1$-skeleton of some combinatorial sphere.

Most of the estimates in this paper will be based on {\it Euler's Formula\/}
for a combinatorial sphere:
if $\Phi$ is a combinatorial sphere, then
$$
\|\Phi(0)\|-\|\Phi(1)\|+\|\Phi(2)\|=2.
$$

Using Euler's Formula, it is easy to prove the following well-known

\proclaim{Proposition 2.1}
For any planar graph\/ $\Gamma$
without loops and without multiple edges$,$
the number of edges of\/ $\Gamma$ is less than
three times the number of its vertices\/$:$
$$
\|\Gamma(1)\|<3\|\Gamma(0)\|.
$$
\endproclaim

\demo{Proof\rm}
Without loss of generality, $\Gamma$ may be assumed to be connected.
(It is enough to prove the above inequality for every connected component of $\Gamma$.)
If $\|\Gamma(0)\|\le2$, then the statement is obvious.
If $\|\Gamma(0)\|\ge3$, then $\Gamma$ is the $1$-skeleton
of some combinatorial sphere $\Phi$,
and the degree of every face of $\Phi$ is at least~$3$.
On one hand,
$$
2\|\Phi(1)\|=\sum_{f\in\Phi(2)}d(f)\ge3\|\Phi(2)\|.
$$
On the other hand, by Euler's Formula,
$$
\|\Phi(0)\|-\|\Phi(1)\|+\|\Phi(2)\|=2.
$$
Therefore,
$$
\aligned
3\|\Phi(0)\|-3\|\Phi(1)\|+3\|\Phi(2)\|+2\|\Phi(1)\|&\ge6+3\|\Phi(2)\|,\\
3\|\Phi(0)\|&\ge\|\Phi(1)\|+6,\\
3\|\Gamma(0)\|&\ge\|\Gamma(1)\|+6>\|\Gamma(1)\|.
\hbox to 0pt{\qed\hss}
\endaligned
$$
\enddemo

An {\it orientation\/} of an edge $e={(} E,\alpha{)}$ of a complex $\Gamma$
is a total order on the two-element set of vertices of~$E$.
There are two possible orientations of every edge,
they are {\it opposite\/} to each other.
An {\it oriented edge\/} is an edge together with one of its orientations.
The set of all oriented edges of a complex $\Gamma$ will be denoted
by~$\hat\Gamma(1)$.
Every morphism $\phi$ of a combinatorial $n$-complex $\Gamma_1\kl$
to a combinatorial $n$-complex $\Gamma_2\kl$, $n\in\{1,2\}$,
naturally defines a function $\hat\Gamma_1\kl(1)\to\hat\Gamma_2\kl(1)$;
denote this function by~$\hat\phi^1$.
Two oriented edges obtained from the same edge by picking
opposite orientations are {\it inverse\/} to each other.
The oriented edge inverse to an oriented edge $e$ is denoted by~$e^{-1}$.
If $e$ is an oriented edge, ${(} E,\alpha{)}$ is a corresponding
non-oriented edge, $E(0)=\{v_1\kl,v_2\kl\}$, and $v_1\kl$ precedes $v_2\kl$ with
respect to the chosen total order on $E(0)$,
then $\bar\alpha^0(v_1\kl)$ is called the {\it tail\/} of $e$ and
$\bar\alpha^0(v_2\kl)$ is called the {\it head\/} of~$e$.
An oriented edge {\it leaves\/} its tail and {\it enters\/} its head.
Clearly, the head of an oriented edge is the tail
of the inverse oriented edge, and vice versa.
For any vertex $v$ in a complex $\Gamma$, the number of oriented edges
leaving $v$ equals the number of oriented edges entering $v$
and equals the degree of~$v$.

Here follow several definitions related to paths
in combinatorial complexes.

A {\it path\/} is a finite sequence of alternating vertices and
oriented edges such that the following conditions hold:
it starts with a vertex and ends with a vertex;
the vertex immediately preceding an oriented edge is its tail;
the vertex immediately following an oriented edge is its head.
The {\it initial vertex\/} of a path is its first vertex,
the {\it terminal vertex\/} of a path is its last vertex,
and an {\it end-vertex\/} of a path is either
its initial or its terminal vertex.
A path {\it starts\/} at its initial vertex and {\it ends\/}
at its terminal vertex.
The {\it length\/} of a path
$p={(} v_0\kl,e_1\kl,v_1\kl,\,\ldots,e_n\kl,v_n\kl{)}$ is $n$;
it is denoted by~$|p|$.
The vertices $v_1\kl,\,\ldots,v_{n-1}\kl$ of this path are called
its {\it intermediate\/} vertices.
A {\it trivial path\/} is a path of length zero.
By abuse of notation, a path of the form ${(} v_1\kl,e,v_2\kl{)}$,
where $e$ is an oriented edge from $v_1\kl$ to $v_2\kl$, shall be denoted by $e$,
and a trivial path ${(} v{)}$ shall be denoted by~$v$.

The {\it inverse path\/} to a path $p$ is defined naturally and
is denoted by~$p^{-1}$.
If the terminal vertex of a path $p_1\kl$ coincides with
the initial vertex of a path $p_2\kl$, then the {\it product\/}
$p_1\kl p_2\kl$ is defined naturally.
A path $s$ is called an {\it initial subpath\/} of a path $p$ if
$p=sq$ for some path~$q$.
A path $s$ is called a {\it terminal subpath\/} of a path $p$ if
$p=qs$ for some path~$q$.

A {\it cyclic path\/} is a path such that its terminal vertex
coincides with its initial vertex.
A {\it cycle\/} is the set of all {\it cyclic shifts\/} of some cyclic path.
The {\it length\/} of a cycle $c$, denoted by $|c|$,
is the length of an arbitrary representative of~$c$.
A {\it trivial cycle\/} is a cycle of length zero.
A path $p$ is a {\it subpath of a cycle\/} $c$ if
for some representative $r$ of $c$ and for some natural $n$,
$p$ is a subpath of $r^n$ (\ie, of the product of $n$ copies of~$r$).
Let $c$ be a cycle in which no oriented edge occurs more than once.
Then a set of paths $S$ is said to {\it cover\/} $c$ if
all the elements of $S$ are nontrivial subpaths of $c$, and every oriented edge
that occurs in $c$ also occurs in some path from~$S$.

A path is {\it reduced\/} if
it does not have a
subpath of the form $ee^{-1}$ where $e$ is an oriented edge.
A cyclic path is {\it cyclically reduced\/} if
it is reduced and its first
oriented edge is not inverse to its last oriented edge.
(For example, all trivial paths are cyclically reduced.)
A cycle is {\it reduced\/} if
it consists of cyclically reduced cyclic paths.
A path is {\it simple\/} if
it is nontrivial, reduced, and
none of its intermediate vertices appears in it more than once.
A cycle is {\it simple\/} if
it consists of simple cyclic paths.

All paths in any graph are naturally partially ordered by
the relation ``is a subpath of.''
A path is called {\it maximal\/} in some set of paths if
it is maximal in this set with respect to the above order
(\ie, if
this path is not a proper subpath of any other path in this set).

Every morphism $\phi$ of combinatorial complexes naturally defines
a function from the set of paths (respectively cycles)
of the first complex to the set of paths (respectively cycles)
of the second.
The image of a given path or cycle under this function
will be referred to as the {\it image relative to\/} $\phi$,
or {\it$\phi$-image}.

An {\it oriented arc\/} in a complex $\Gamma$ is a simple path
all intermediate vertices of which have degree $2$ in~$\Gamma$.
A {\it non-oriented arc}, or simply {\it arc}, may be defined as a pair
of mutually inverse oriented arcs.
The {\it length\/} of a non-oriented arc $u$ is defined as the length
of either of the two associated oriented arcs and is denoted by~$|u|$.
The concepts of a {\it subarc}, a {\it maximal arc},
and so on, are self-explanatory.
Sometimes edges will be viewed as arcs, and oriented edges as oriented arcs.
The {\it set of edges\/} of an oriented arc (respectively of an arc) is the set of all
edges that with some orientation occur in that oriented arc
(respectively in either of the two corresponding oriented arcs).
An {\it intermediate vertex\/} of an arc is an intermediate vertex
of either of the corresponding oriented arcs.
Two arcs, or oriented arcs, are said to {\it overlap\/} if
their sets of edges are not disjoint.
A set of arcs, or oriented arcs, is called {\it non-overlapping\/} if
any two distinct elements of this set have disjoint sets of edges.
An arc (or an edge) $u$ {\it lies\/} on a path $p$ if
at least one of the oriented arcs (respectively oriented edges)
associated with $u$ is a subpath of~$p$.
A set of arcs $A$ {\it covers\/} a set of arcs (or edges) $B$ if
every element of $B$ is a subarc of some element of~$A$.

Combinatorial complexes representing orientable surfaces
may be oriented in a combinatorial manner.
To define the combinatorial notion of {\it orientation},
consider first an arbitrary combinatorial circle~$\Gamma$.
Consider a function that for every edge $e$ of $\Gamma$ chooses
an orientation of~$e$.
Call an edge with a chosen orientation a {\it chosen oriented edge}.
Say that such a function chooses {\it coherent orientations\/} if
every vertex of $\Gamma$ is the tail of exactly one
of the chosen oriented edges and the head of exactly one
of the chosen oriented edges.
An {\it orientation\/} of a combinatorial circle is a function
on the set of its edges that chooses coherent orientations of the edges.
An {\it oriented combinatorial circle\/} is a combinatorial circle
together with one of its two possible orientations.
Note that in an oriented combinatorial circle $\Gamma$
there is a unique simple cycle whose oriented edges are the chosen ones.
Call this cycle the {\it chosen cycle\/} of~$\Gamma$.

Now, consider a face $f={(} F,\beta{)}$ of a $2$-complex~$\Phi$.
An {\it orientation\/} of $f$ is an orientation of the combinatorial circle~$F$.
An {\it oriented face\/} is a face together with one of its orientations.
The set of all oriented faces of a complex $\Phi$ will be denoted
by~$\hat\Phi(2)$.
(Note that $\|\hat\Phi(2)\|=2\|\Phi(2)\|$.)
Every morphism $\phi$ of a combinatorial $2$-complex $\Phi_1\kl$
to a combinatorial $2$-complex $\Phi_2\kl$
naturally defines a function $\hat\Phi_1\kl(2)\to\hat\Phi_2\kl(2)$;
denote this function by~$\hat\phi^2$.
Two oriented faces obtained from the same face by picking
{\it opposite\/} orientations are {\it inverse\/} to each other.
The oriented face inverse to an oriented face $f$ is denoted by~$f^{-1}$.
Each oriented face has a uniquely defined {\it boundary cycle\/}:
if $f$ is an oriented face of $\Phi$ and ${(} F,\beta{)}$ is its underlying
non-oriented face, then the boundary cycle of $f$ is the cycle in $\Phi$
which is the $\beta$-image of the chosen cycle of~$F$.
Boundary cycles of mutually inverse oriented faces are
mutually inverse; they are also called {\it boundary cycles\/} of the
corresponding non-oriented face.

Consider now a combinatorial surface~$\Phi$.
Consider a function $\theta$ that for every face $f$ of $\Phi$ chooses
an orientation of~$f$.
Call a face with a chosen orientation a {\it chosen oriented face}.
For every oriented edge $e$ of $\Phi$, count the number of all such
ordered pairs ${(} f,x{)}$ that $f={(} F,\beta{)}$ is a face of $\Phi$,
$x$ is an oriented edge of $F$, $x$ is chosen with respect
to the orientation $\theta(f)$ of $F$, and $e=\hat\beta^1(x)$
(in particular, $f$ is incident to~$e$).
The function $\theta$ is said to choose {\it coherent orientations\/} if
for every oriented edge of $\Phi$ the number defined above
is either $1$ or~$0$.
An {\it orientation\/} of a combinatorial surface is a choice of
coherent orientations of all of its faces.
A combinatorial surface that admits orientation is called {\it orientable}.
The fact of the matter is that a combinatorial surface is orientable if and only if
it represents an orientable surface,
but this fact is not used in the paper.
There exist exactly two orientations of any connected orientable combinatorial surface.
The {\it contour\/} of a (non-oriented) face in an oriented combinatorial surface
is the boundary cycle of the corresponding chosen oriented face.
An {\it isomorphism\/} of oriented combinatorial surfaces $A$ and $B$
is an isomorphism of $A$ and $B$ as $2$-complexes which
additionally preserves the chosen orientations of the faces.

The only $2$-complexes used in this paper are subcomplexes
of combinatorial spheres.
Connected components of such subcomplexes
with some additional structure will be called {\it maps}.

\definition{Definition}
A {\it spherical map\/} is an oriented combinatorial sphere.
\enddefinition

Consider an oriented combinatorial sphere~$\Delta$.
By removing some of its faces, one can obtain a subcomplex $\Phi$
of the underlying $2$-complex of~$\Delta$.
Note that such a subcomplex $\Phi$ is connected
($\Phi$ and $\Delta$ have the same $1$-skeleton).
If the contours of the removed faces and
the orientations of the faces of $\Phi$ induced from $\Delta$ are known,
then the oriented combinatorial sphere $\Delta$ may be reconstructed from $\Phi$
up to isomorphism.
(Just attach new faces to cover the ``holes'' in $\Phi$
and orient all faces appropriately.)

\definition{Definition}
A {\it nontrivial map\/} $\Phi$ is any
$2$-complex that can be obtained from some
oriented combinatorial sphere $\Delta$ in the way described above,
together with the inherited orientation of its faces
and with the set of the contours of the removed faces of $\Delta$,
called the {\it contours\/} of~$\Phi$.
Such a spherical map $\Delta$ is called a {\it spherical closure\/}
of the given nontrivial map~$\Phi$.
If $\Phi$ is a nontrivial map and $\Delta$ is its spherical closure,
then those faces of $\Delta$ that are in the same time faces of $\Phi$
are called {\it proper}, and the other ones are called {\it improper}.
A {\it trivial map\/} is a $2$-complex consisting
of a single vertex together with the one-element set consisting of
the only (trivial) cycle in this $2$-complex.
This (trivial) cycle is called the {\it contour\/} of this trivial map.
A map is called {\it degenerate\/} if
its set of faces is empty.
\enddefinition

\remark{Remark\/ \rm2.1}
At this point, one may want to redefine a spherical map as
an oriented combinatorial sphere together with the empty set
(of its contours).
\endremark

The following definition is used to distinguish ``pathological'' maps.

\definition{Definition}
An {\it exceptional map\/} is a spherical map,
the $1$-skeleton of whose underlying $2$-complex is a combinatorial circle
(such a map must consist of two faces ``glued'' together along their boundary cycles).
\enddefinition

In a non-exceptional spherical map, any two distinct maximal arcs
do not overlap.

Any connected subcomplex of any map naturally inherits
a map structure (even if it was obtained by removing not only
faces, but also edges and vertices).

\definition{Definition}
A {\it submap\/} of a map is a connected subcomplex with
the inherited map structure.
\enddefinition

Note: in general, a spherical map $\Delta$ is not a spherical closure
of a given nontrivial submap of~$\Delta$.

In every map, each oriented edge appears either in the contour
of some face of the map or in a contour of the map.

\definition{Definition}
A map is {\it simple\/} if
all of its contours are simple cycles without common vertices.
\enddefinition

\remark{Remark\/ \rm2.2}
If a map $\Delta$ is simple,
then every edge of $\Delta$ is incident with some face of~$\Delta$.
\endremark

\definition{Definition}
A {\it disc map\/} is a map with exactly one contour.
\enddefinition

\remark{Remark\/ \rm2.3}
A disc map is simple if and only if its underlying $2$-complex is a combinatorial disc.
\endremark

The contours of faces in a disc map may be thought of as
{\it oriented counterclockwise}, and the contour of the map itself
may be thought of as
{\it oriented clockwise\/}
(this convention corresponds to the way a disc map is usually pictured).
Divide the set of all simple cycles in a disc map into
{\it oriented clockwise\/} and {\it oriented counterclockwise\/} accordingly.

An oriented arc $u$ of a map $\Phi$
is {\it incident\/} to a face $f$ of $\Phi$ if
it is a subpath of one of the boundary cycles of~$f$.
An arc is {\it incident\/} to $f$ if
one/both of the corresponding oriented arcs are incident to~$f$.
An arc incident to a face $f$ is also called an {\it arc\/} of~$f$.
An arc which is incident to some face is either {\it internal\/}
(does not lie on any of the contours of the map) or {\it external\/}
(lies on some contour of the map).
Internal oriented arcs divide into {\it inter-facial\/}
(incident to two different faces)
and {\it intro-facial\/} (incident to only one face).
An internal arc $u$ is {\it between\/} faces $f_1\kl$ and $f_2\kl$ if
$\{f_1\kl,f_2\kl\}$ is the set of all faces incident to~$u$.
Every edge may be considered as an arc.
Therefore, it makes sense to say that an edge incident to some face
is internal, external, inter-facial, or intro-facial.
An oriented arc or an oriented edge is internal, external,
inter-facial or intro-facial if
the corresponding non-oriented arc or edge is such.
Every intro-facial oriented arc $u$ in a disc map is either {\it outward\/}
($u$ is an initial subpath of some simple path with the terminal vertex on the contour of the map),
or {\it inward\/} ($u^{-1}$ is outward).

\definition{Definition}
Let $\Delta$ be a map, $c$ be a nontrivial cycle in $\Delta$,
$\Delta'$ be a simple disc map.
Say that $c$ {\it cuts\/} $\Delta'$ out of $\Delta$ if
there exists a morphism $\zeta$ of the underlying $2$-complex of $\Delta'$
to the underlying $2$-complex of $\Delta$
preserving the chosen orientations of the faces
such that $\bar\zeta^2\:\Delta'(2)\to\Delta(2)$ is injective
and $c$ is the $\zeta$-image of the contour of~$\Delta'$.
Call such a morphism $\zeta$ a {\it pasting morphism}.
\enddefinition

\remark{Remark\/ \rm2.4}
If a cycle $c$ cuts a simple disc map $\Delta'$ out of a map $\Delta$,
then $\Delta'$ is essentially determined by $\Delta$ and~$c$.
\endremark

It is also common to call a simple disc map that is cut out of a given map $\Delta$
just a ``simple disc submap'' of $\Delta$,
but it is not compatible with the definition of a submap given in this paper.

If $\Delta_0\kl$ is a simple disc submap of a map $\Delta$,
then the contour of $\Delta_0\kl$ cuts $\Delta_0\kl$ out of $\Delta$,
and the natural morphism of $\Delta_0\kl$ to $\Delta$
is the corresponding pasting morphism.
In a disc map, any simple cycle oriented clockwise is the contour of
some simple disc submap; therefore, it cuts out that simple disc submap.

Let $\Delta$ be a combinatorial surface
or a subcomplex of a combinatorial surface.
Call two distinct faces of $\Delta$ {\it contiguous\/} if
there exist an edge incident with both of them.
(The notion of contiguity will only be used for distinct faces.)
Let $F$ be a nonempty set of faces of~$\Delta$.
Let $S$ be the set of all (non-ordered) pairs of distinct
contiguous faces in the set~$F$.
Clearly, there exist a graph $\Gamma$ and bijections
$$
\alpha_0\kl\:F\to\Gamma(0),\quad
\alpha_1\kl\:S\to\Gamma(1),
$$
such that a face $f\in F$ belongs to a pair $P\in S$ if and only if
the vertex $\alpha_0\kl(f)$ is incident to the edge~$\alpha_1\kl(P)$.
Indeed, choose vertices of $\Gamma$ in bijective correspondence with faces of the set $F$,
and connect any two vertices that correspond to distinct contiguous faces with an edge.

\definition{Definition}
Such a graph $\Gamma$ is called a {\it contiguity graph\/}
for the set $F$ in~$\Delta$.
\enddefinition

\remark{Remark\/ \rm2.5}
A contiguity graph is unique up to graph isomorphism.
Obviously, any contiguity graph cannot have loops or multiple edges.
A contiguity graph of the set of all faces in a combinatorial surface is connected.
It is easy to see that a contiguity graph of any nonempty set of faces
in (any subcomplex of) a combinatorial sphere is a planar graph.
\endremark

Consider a combinatorial sphere~$\Delta$.
Let $E$ be a nonempty non-overlapping set of arcs of $\Delta$
satisfying the following condition:
for any two distinct elements $e_1\kl,e_2\kl$ of $E$, there exists
a finite sequence ${(} x_0\kl,\,\ldots,x_n\kl{)}$ of elements
of $E$ such that $x_0\kl=e_1\kl$, $x_n\kl=e_2\kl$,
and for every $i=0$, \dots, $n-1$, the arcs $x_i\kl$
and $x_{i+1}\kl$ are both incident to some face of~$\Delta$.
Let $F$ be the set of all faces of $\Delta$ that are
incident to at least one element of~$E$.
Let $\Psi$ be the $2$-complex obtained from $\Delta$
by removing all the faces that belong to $F$, and
all the edges and all the intermediate vertices
of all the arcs that belong to~$E$.
Note that all the end-vertices of all the arcs from $E$ are vertices of~$\Psi$.
Let $V$ be the set of all connected components of~$\Psi$.
Let $\hat E$ be the set of all oriented arcs associated with arcs from~$E$.

\definition{Definition}
An {\it estimating\/ $2$-complex\/} for the set $E$ in the combinatorial sphere $\Delta$
is a combinatorial sphere $\Phi$ such that there exist bijections
$$
\alpha_0\Kl\:V\to\Phi(0),\quad
\alpha_1'\:\hat E\to\hat\Phi(1),\quad
\alpha_2\Kl\:F\to\Phi(2),
$$
satisfying the following properties:
\medskip

\item{$\bullet$}
The bijection $\alpha_1'$ takes mutually inverse oriented arcs of $\Delta$ to
mutually inverse oriented edges of~$\Phi$.

\item{$\bullet$}
For each $v\in V$, $e\in\hat E$,
the vertex $\alpha_0\kl(v)$ is the tail (respectively head)
of the oriented edge $\alpha_1'(e)$ if and only if
the initial vertex (respectively terminal vertex)
of the oriented arc $e$ belongs to the component~$v$.

\item{$\bullet$}
If one of the boundary cycles of $f\in F$ has a representative
of the form
$$
\displaywidth=\hsize
\displayindent=0pt
x_0\kl y_1\kl x_1\kl\cdots y_n\kl x_n\kl,
$$
where each path $x_i\kl$ is wholly contained in
some $v_i\kl\in V$, and $y_1\kl,\,\ldots,y_n\kl\in\hat E$, then
$$
\bigl{(}\alpha_0\Kl(v_0\Kl),\alpha_1'(y_1\Kl),\alpha_0\Kl(v_1\Kl),\,\ldots,
\alpha_1'(y_n\Kl),\alpha_0\Kl(v_n\Kl)\bigr{)}
$$
is a representative of a boundary cycle of the face~$\alpha_2\kl(f)$.
\enddefinition

\remark{Remark\/ \rm2.6}
It may be shown under the above assumptions that an estimating $2$-complex exists
and is unique up to isomorphism.
\endremark

Informally speaking, an estimating $2$-complex for the set $E$ in $\Delta$
may be obtained as follows.
First, take all the edges that do not belong to elements of $E$
and all the faces that are not incident to elements of $E$ and collapse them into points.
These points will be some of the vertices of the estimating $2$-complex.
Second, erase all the intermediate vertices from the elements of $E$
so that every arc from $E$ becomes an edge.
The condition imposed on $E$ ensures that the obtained $2$-complex
is a combinatorial sphere.

\topinsert
\vskip .125in
\centerline{\epsfbox{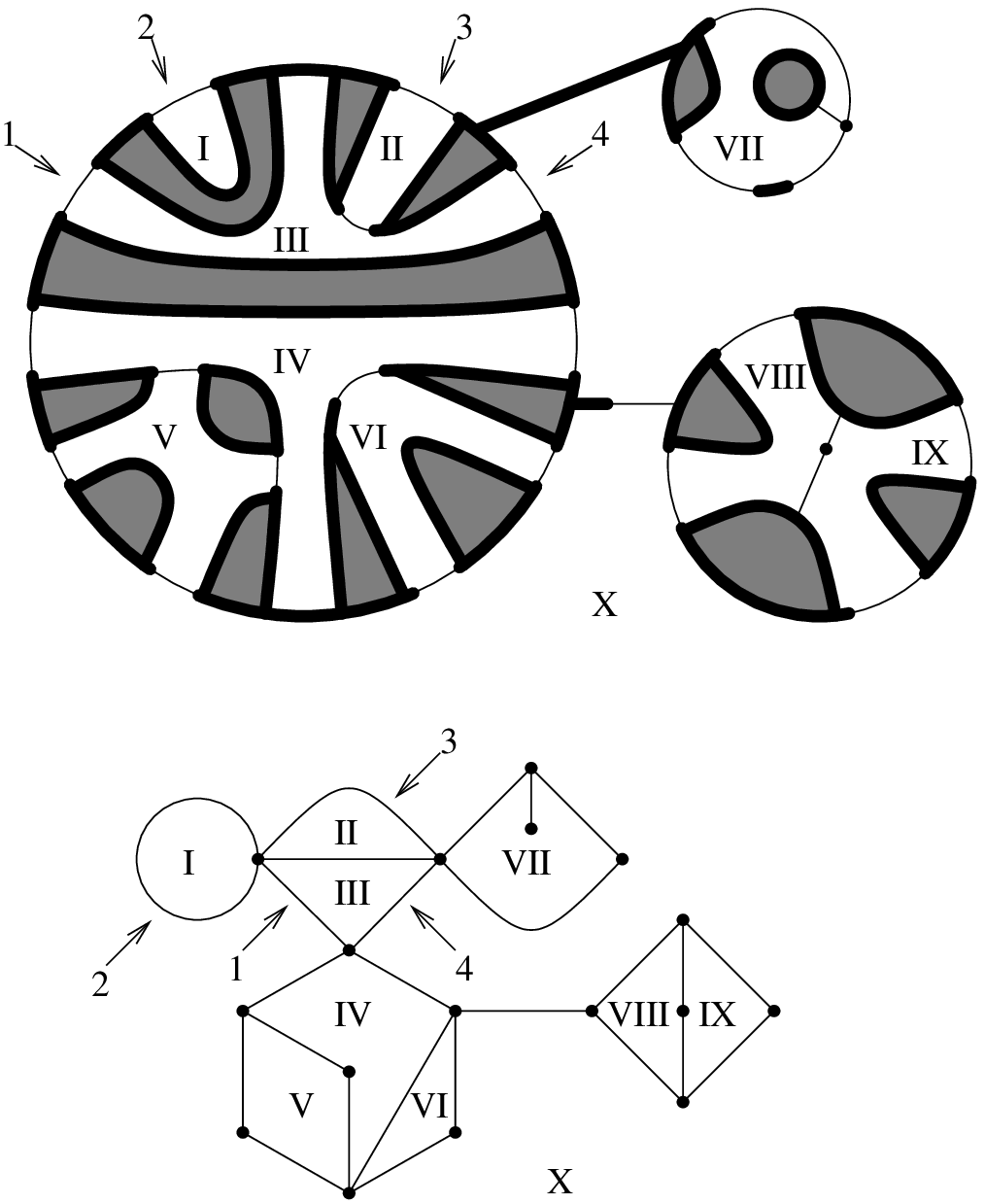}}
\botcaption{Fig.~1: Combinatorial sphere and estimating $2$-complex.}
\endcaption
\endinsert

Figure~1 shows an example of a combinatorial sphere $\Delta$ and an estimating
$2$-complex $\Phi$ for some non-overlapping set of arcs $E$ in~$\Delta$.
The upper part of the picture shows a stereographic projection of $\Delta$,
in which the elements of $E$ are indicated by the thinner lines.
The lower part shows a stereographic projection of~$\Phi$.
The roman numerals indicate the correspondence between the faces of $\Delta$
incident with the arcs from $E$ ($10$ faces) and the faces of~$\Phi$.
The arabic numerals indicate the correspondence between some arcs from $E$
and the corresponding edges of~$\Phi$
(the total number of elements in $E$ is~$25$).
The parts of $\Delta$ that are to be ``collapsed'' are shown as shaded regions and
thicker lines (there are $15$ such ``parts'').
There are also $2$ vertices in $\Phi$ that correspond to vertices of~$\Delta$.

The estimations used in the proof of the main theorem
are based on Euler's Formula applied to estimating $2$-complexes.
Objects similar to estimating $2$-complexes were extensively
used in~\cite{O};
there they were called {\it estimating graphs}.
\par

\head
3. Maps with selection
\endhead

The proofs of Theorems 1 and 2 given in the end of this paper
are based on constructing
group presentations such that some sets of subwords
of the defining relations satisfy certain small-cancellation-type conditions,
though the presentations themselves do not need to satisfy classical
small cancellation conditions.
This is why it is convenient to consider maps and
van Kampen diagrams together with some systems of selected subpaths of
the contours of their faces.

\definition{Definition}
A {\it selection\/} on a face $f$ of a map $\Delta$ is a set of nontrivial
reduced subpaths of the contour of $f$
such that for each path in this set, all of its nontrivial subpaths
belong to the set as well (\ie, the set is closed under taking nontrivial subpaths).
A {\it selection\/} on a map $\Delta$ is a set of nontrivial reduced subpaths
of the contours of faces of $\Delta$ which is closed under taking nontrivial subpaths.
An {\it S-map}, or {\it map with selection}, is a map together with a selection on it.
\enddefinition

A path in an S-map is {\it selected\/} if
it belongs to the selection.
An oriented edge in an S-map is {\it selected\/} if
it is selected as a path.
A path or an oriented edge is {\it double-selected\/} if
it is selected along with its inverse.
An edge or an arc is {\it double-selected\/} if
it is internal and both of the corresponding oriented edges or arcs are selected.
An external edge is {\it selected\/} if
one of the associated oriented edges is selected.

If $\zeta$ is a morphism of the underlying $2$-complex of a map $\Delta_1\kl$
to the underlying $2$-complex of a map $\Delta_2\kl$,
and $\zeta$ preserves the chosen orientations of the faces,
then any selection on $\Delta_2\kl$
naturally induces a selection on $\Delta_1\kl$ via $\zeta$:
a path $p$ in $\Delta_1\kl$ is selected if and only if
its $\zeta$-image in $\Delta_2\kl$ is selected
(note that if the $\zeta$-image of a path is reduced,
then the path itself is reduced).
Say that a cycle $c$ {\it cuts\/} a simple disc S-map $\Delta'$
out of an S-map $\Delta$ if
$c$ cuts the underlying map of $\Delta'$
out of the underlying map of $\Delta$,
and the selection on $\Delta'$ is the one induced from $\Delta$
via some pasting morphism.
Only those pasting morphisms of the underlying map of $\Delta'$
that induce the given selection on $\Delta'$
are called {\it pasting morphisms\/} of the S-map $\Delta'$ to~$\Delta$.

Any submap of an S-map has a naturally induced selection.
Say that an S-map $\Delta_0\kl$ is an {\it S-submap\/} of an S-map $\Delta$ if
$\Delta_0\kl$ is a submap of $\Delta$ together with
the induced selection.

By abuse of notation, S-maps sometimes will be called just maps.

\definition{Definition}
An S-map $\Delta$ is said to satisfy the
{\it condition $\Cal Z(n)$}, $n\in\Bbb N\cup\{0\}$, if
for every set $S$ of selected paths that covers the contour of
some simple disc submap of $\Delta$,
the set $S$ has more than $n$ maximal elements;
in particular, $S$ has at least $n+1$ element.
\enddefinition

\proclaim{Proposition 3.1}
Let\/ $\Phi$ be a non-degenerate disc map$.$
Suppose a set\/ $S$ consists of\/ $n$ reduced paths and
covers the contour of\/~$\Phi.$
Then there exists a maximal simple disc submap\/ $\Phi_1\kl$ of\/ $\Phi$
whose contour may be covered by\/ $n$ or fewer subpaths
of the paths comprising\/~$S.$
\endproclaim

\demo{Proof\rm}
Let $c$ be the result of cyclic reduction of the contour of~$\Phi$.
Since $\Phi$ is a non-degenerate disc map,
the cycle $c$ is the contour of a non-degenerate disc map~$\Phi'$;
in particular, $c$ is nontrivial.
Let $q_1\kl$, \dots, $q_{n'}\kl$ be such subpaths of $c$ that
the product $q_1\kl\cdots q_{n'}\kl$ is a cyclic path representing $c$,
each $q_i\kl$ is a nontrivial subpath of some element of $S$,
and each element of $S$ has at most one of the paths $q_1\kl$, \dots, $q_{n'}\kl$
as a subpath.
It is not hard to prove that such paths $q_1\kl$, \dots, $q_{n'}\kl$ exist.
Note that $n'\le n$.

If $\Phi'$ is simple, then $\Phi_1\kl=\Phi'$ is a desired submap.
Consider the case when $\Phi'$ is not simple.
Since the map $\Phi'$ is non-degenerate, it has two distinct
maximal simple disc submaps $\Phi_1\kl$ and $\Phi_2\kl$,
whose contours have representatives which are subpaths of
the contour of~$\Phi'$.
Let $p_i\kl$ be the representative of the contour of $\Phi_i\kl$
which is a subpath of the contour of $\Phi'$, $i=1$,~$2$.
Let $v_i\kl$ be the initial vertex of $q_i\kl$, $i=1$, \dots,~$n'$.
Let $w_i\kl$ be the initial (and the terminal) vertex of $p_i\kl$, $i=1$,~$2$.
Let $k_i\kl$ be the number of vertices from the set
$\{v_1\kl,\,\ldots,v_{n'}\kl,w_1\kl,w_2\kl\}$
which occur in the contour of $\Phi_i\kl$, $i=1$,~$2$.
It is easy to see that either $k_1\kl\le n'$ or $k_2\kl\le n'$.
For the sake of definiteness, assume that $k_1\kl\le n'$.
Then the contour of $\Phi_1\kl$ has a representative of the form
$q_1^{(1)}\cdots q_{k_1\kl}^{(1)}$ where
every $q_i^{(1)}$ is a nontrivial subpath of some~$q_j\kl$.
The set $\{q_1^{(1)},\,\ldots,q_{k_1\kl}^{(1)}\}$ covers the contour of~$\Phi_1\kl$.
\qed
\enddemo

\proclaim{Corollary}
Let\/ $\Delta$ be an S-map$.$
Suppose the contour of some non-degenerate
disc submap of\/ $\Delta$ is covered by
a set of\/ $n$ or fewer selected paths$.$
Then\/ $\Delta$ does not satisfy the condition\/ $\Cal Z(n).$
\endproclaim

\head
4. Condition $\Cal A$
\endhead

Roughly speaking, the condition $\Cal A$ defined in this section
is a generalization of the small cancellation condition $C'$
formulated in terms of the underlying maps of reduced van Kampen diagrams,
rather than in terms of defining relations.
(About the condition $C'$ and van Kampen diagrams, see Section 11.)

Recall that the degree $d(f)$ of a face $f$ of a map is defined as
the length of the contour of~$f$.

\definition{Definition}
Let $\Delta$ be an S-map.
Let $\bar a={(} k_1\kl,k_2\kl,k_3\kl;\lambda_1\kl,\lambda_2\kl,\lambda_3\kl,\lambda_4\kl{)}$
where
$$
k_1\kl,k_2\kl,k_3\kl\:\Delta(2)\to\Bbb N\cup\{0\},\quad
\lambda_1\kl,\lambda_2\kl,\lambda_3\kl,\lambda_4\kl\:\Delta(2)\to[0,1].
$$
The S-map $\Delta$ is said to satisfy the {\it condition\/}
$\Cal A(\bar a)$ if
it satisfies the following seven conditions:
\medskip

\par\hang\noindent{$\Cal A_1\kl(k_1\kl)$}
For each face $f$ of $\Delta$,
there exists at least one selected subpath of the contour of $f$,
and the number of maximal selected subpaths
of the contour of $f$ does not exceed~$k_1\kl(f)$
(note that if all nontrivial subpaths are selected,
then there is no maximal selected subpath).

\par\hang\noindent{$\Cal A_2\kl(k_2\kl)$}
For each face $f$ of $\Delta$,
there are at most $k_2\kl(f)$
subpaths of the contour of $f$ that
are of the form $p_1\kl qp_2\kl$ where $p_1\kl$ and $p_2\kl$ are
maximal selected paths, $q$ is a nontrivial path,
and no subpath of $q$ is selected.

\par\hang\noindent{$\Cal A_3\kl(k_3\kl)$}
For each face $f$ of $\Delta$,
there are at most $k_3\kl(f)$
subpaths of the contour of $f$ that
are of the form $p_1\kl p_2\kl$ where $p_1\kl$ and $p_2\kl$ are
maximal selected paths, and $p_1\kl p_2\kl$ is not reduced.

\par\hang\noindent{$\Cal A_4\kl(\lambda_1\kl)$}
For each face $f$ of $\Delta$,
if $S$ is the number of non-selected oriented
edges in the contour of $f$, then
$$
\displaywidth=\hsize
\displayindent=0pt
S\le\lambda_1\kl(f)d(f).
$$

\par\hang\noindent{$\Cal A_5\kl(\lambda_2\kl)$}
For each face $f$ of $\Delta$,
if $u$ is a double-selected intro-facial arc incident to $f$,
then
$$
\displaywidth=\hsize
\displayindent=0pt
|u|\le\lambda_2\kl(f)d(f).
$$

\par\hang\noindent{$\Cal A_6\kl(\lambda_2\kl,\lambda_3\kl)$}
For every two distinct faces $f_1\kl$, $f_2\kl$ of $\Delta$,
if $U$ is a non-overlapping set of double-selected arcs
between $f_1\kl$ and $f_2\kl$
which covers the set of all double-selected edges
between $f_1\kl$ and $f_2\kl$,
then
$$
\displaywidth=\hsize
\displayindent=0pt
\sum_{u\in U}|u|
\le\|U\|\cdot\min\bigl\{\lambda_2\kl(f_1\kl)d(f_1\kl),\lambda_2\kl(f_2\kl)d(f_2\kl)\bigr\}
+\min\bigl\{\lambda_3\kl(f_1\kl)d(f_1\kl),\lambda_3\kl(f_2\kl)d(f_2\kl)\bigr\}.
$$

\par\hang\noindent{$\Cal A_7\kl(\lambda_2\kl,\lambda_4\kl)$}
For each face $f$ of $\Delta$,
if $p$ is a simple selected subpath of the contour of a distinct face of $\Delta$,
$U$ is a non-overlapping set of double-selected arcs
incident to $f$ and lying on $p$,
and $U$ covers the set of all double-selected edges
that are incident to $f$ and lie on $p$,
then
$$
\displaywidth=\hsize
\displayindent=0pt
\sum_{u\in U}|u|
\le\bigl(\|U\|\lambda_2\kl(f)+\lambda_4\kl(f)\bigr)d(f).
$$
In particular, if a double-selected inter-facial arc $u$ is incident to a face $f$,
then
$$
|u|\le\bigl(\lambda_2\kl(f)+\lambda_4\kl(f)\bigr)d(f).
$$

\medskip

\noindent
A map $\Delta$ is said to satisfy the condition $\Cal A(\bar a)$ if
there exists a selection on $\Delta$ such that $\Delta$ with this selection
satisfies $\Cal A(\bar a)$.
\enddefinition

\remark{Remark\/ \rm4.1}
If for each face $f$ of $\Delta$, the length of every double-selected arc
incident to $f$ is at most $\lambda_2\kl(f)d(f)$,
then $\Delta$ automatically satisfies the conditions
$\Cal A_5\kl(\lambda_2\kl)$, $\Cal A_6\kl(\lambda_2\kl,\lambda_3\kl)$, $\Cal A_7\kl(\lambda_2\kl,\lambda_4\kl)$.
If $\Delta$ satisfies $\Cal A_1\kl(k_1\kl)$,
then it also satisfies $\Cal A_2\kl(k_1\kl)$ and $\Cal A_3\kl(k_1\kl)$.
\endremark

If $\Delta_1\kl$ and $\Delta_2\kl$ are two maps,
$\zeta$ is a morphism from $\Delta_1\kl$ to $\Delta_2\kl$ which preserves
the chosen orientations of the faces,
and
$$
\gathered
k_1\kl,k_2\kl,k_3\kl\:\Delta_2\kl(2)\to\Bbb N\cup\{0\},\quad
\lambda_1\kl,\lambda_2\kl,\lambda_3\kl,\lambda_4\kl\:\Delta_2\kl(2)\to[0,1],\\
\bar a={(} k_1\kl,k_2\kl,k_3\kl;\lambda_1\kl,\lambda_2\kl,\lambda_3\kl,\lambda_4\kl{)},
\endgathered
$$
then define
$$
\gathered
k_1\kl.\zeta=k_1\kl\circ\bar\zeta^2,\quad
k_2\kl.\zeta=k_2\kl\circ\bar\zeta^2,\quad
k_3\kl.\zeta=k_3\kl\circ\bar\zeta^2,\\
\lambda_1\kl.\zeta=\lambda_1\kl\circ\bar\zeta^2,\quad
\lambda_2\kl.\zeta=\lambda_2\kl\circ\bar\zeta^2,\quad
\lambda_3\kl.\zeta=\lambda_3\kl\circ\bar\zeta^2,\quad
\lambda_4\kl.\zeta=\lambda_4\kl\circ\bar\zeta^2,\\
\bar a.\zeta={(} k_1\kl.\zeta,k_2\kl.\zeta,k_3\kl.\zeta;
\lambda_1\kl.\zeta,\lambda_2\kl.\zeta,\lambda_3\kl.\zeta,\lambda_4\kl.\zeta{)},
\endgathered
$$
where $\bar\zeta^2$ is the function $\Delta_1\kl(2)\to\Delta_2\kl(2)$ associated with~$\zeta$.

\remark{Remark\/ \rm4.2}
If an S-map $\Delta$ satisfies the condition $\Cal A(\bar a)$,
a simple disc S-map $\Delta'$ is cut out of $\Delta$,
and $\zeta$ is a pasting morphism,
then the S-map $\Delta'$ satisfies the condition $\Cal A(\bar a.\zeta)$.
A similar statement is true for each of the conditions
$\Cal A_1\kl$--$\Cal A_7\kl$ separately.
\endremark

\head
5. Estimating lemmas
\endhead

This section contains two lemmas used in the proof of the main theorem
for estimating the total number of double-selected internal edges of an S-map.
The proofs of these lemmas use one well-known combinatorial result
published by Philip Hall in 1935~\cite{H}.

If $X$ is a set, then $\eusm P(X)$ shall denote the set of all subsets of $X$
(the power set of~$X$).

\proclaim{Lemma (P\. Hall, 1935)}
Let\/ $A$ and\/ $B$ be two finite sets$.$
Let\/ $f$ be a function from\/ $A$ to\/ $\eusm P(B).$
Let a function\/ $F\:\eusm P(A)\to\eusm P(B)$ be defined by\/
$F(X)=\bigcup_{x\in X}f(x).$
Then the following are equivalent\/$:$
\medskip

\itemitem{$(a)$}
There exists an injection\/ $h\:A\to B$ such that
for all\/ $x\in A,$ $h(x)\in f(x).$
\smallskip

\itemitem{$(b)$}
For each subset\/ $X$ of\/ $A,$ $\|X\|\le\|F(X)\|.$

\endproclaim

The proof of this fact given here appears to be more concise then the original one:

\demo{Proof\rm}
The implication $(a)\Rightarrow(b)$ is obvious.
To prove the converse implication, induct on~$\|A\|$.

If $\|A\|=0$, then $A=\emptyset$ and the conclusion of this lemma
is obvious (take $h=\emptyset$).

Let $n$ be a natural number.
Suppose the implication $(a)\Leftarrow(b)$ holds under the additional
assumption that $\|A\|<n$.
Now, assume that $\|A\|=n$ and suppose that $(b)$ holds.

Case 1:
there exists a proper nonempty subset $A_1\kl$ of the set $A$ such that
$$
\|F(A_1\kl)\|=\|A_1\kl\|.
$$
Then let $B_1\kl=F(A_1\kl)$, $A_2\kl=A\setminus A_1\kl$, $B_2\kl=B\setminus B_1\kl$.
Note that $\|A_1\kl\|<n$ and $\|A_2\kl\|<n$.
Let $f_1\kl=f|_{A_1\kl}$, $F_1\kl=F|_{A_1\kl}$.
Note that $f_1\kl\:A_1\kl\to\eusm P(B_1\kl)$ and
$F_1\kl\:\eusm P(A_1\kl)\to\eusm P(B_1\kl)$.
Let
$f_2\kl\:A_2\kl\to\eusm P(B_2\kl)$,
$F_2\kl\:\eusm P(A_2\kl)\to\eusm P(B_2\kl)$
be defined by
$$
f_2\kl(x)=f(x)\cap B_2\kl,\qquad
F_2\kl(X)=F(X)\cap B_2\kl.
$$
Obviously,
$$
\bigl(\forall X\subset A_1\kl\bigr)\bigl(\|X\|\le\|F_1\kl(X)\|\bigr).
$$
It is also easy to see that
$$
\bigl(\forall X\subset A_2\kl\bigr)\bigl(\|X\|\le\|F_2\kl(X)\|\bigr).
$$
Therefore, by the inductive assumption,
there exist injections $h_1\kl\:A_1\kl\to B_1\kl$ and $h_2\kl\:A_2\kl\to B_2\kl$
such that
$$
\bigl(\forall x\in A_1\kl\bigr)\bigl(h_1\kl(x)\in f_1\kl(x)\bigr)\quad\text{and}\quad
\bigl(\forall x\in A_2\kl\bigr)\bigl(h_2\kl(x)\in f_2\kl(x)\bigr).
$$
Clearly, $h=h_1\kl\cup h_2\kl$ is an injection $A\to B$ such that
$$
\bigl(\forall x\in A\bigr)\bigl(h(x)\in f(x)\bigr).
$$

Case 2:
for each proper nonempty subset $A_1\kl$ of the set $A$,
$$
\|F(A_1\kl)\|\ge\|A_1\kl\|+1.
$$
Then take an arbitrary $x_0\kl\in A$ and an arbitrary $y_0\kl\in f(x_0\kl)$.
Let $A_1\kl=A\setminus\{x_0\kl\}$, $B_1\kl=B\setminus\{y_0\kl\}$.
Note that $\|A_1\kl\|=n-1$.
Let
$f_1\kl\:A_1\kl\to\eusm P(B_1\kl)$,
$F_1\kl\:\eusm P(A_1\kl)\to\eusm P(B_1\kl)$
be defined by
$$
f_1\kl(x)=f(x)\setminus\{y_0\kl\},\qquad
F_1\kl(X)=F(X)\setminus\{y_0\kl\}.
$$
It is obvious that
$$
\bigl(\forall X\subset A_1\kl\bigr)\bigl(\|X\|\le\|F_1\kl(X)\|\bigr).
$$
Therefore, by the inductive assumption,
there exists an injection $h_1\kl\:A_1\kl\to B_1\kl$
such that
$$
\bigl(\forall x\in A_1\kl\bigr)\bigl(h_1\kl(x)\in f_1\kl(x)\bigr).
$$
Let $h\:A\to B$ be the extension of $h_1\kl$ to the whole of $A$ by
putting $h(x_0\kl)=y_0\kl$.
Clearly, $h$ is an injection satisfying
$$
\bigl(\forall x\in A\bigr)\bigl(h(x)\in f(x)\bigr).
$$

Since either Case 1 or Case 2 must take place,
the implication $(a)\Leftarrow(b)$ holds when $\|A\|=n$.
The inductive step is done.
\qed
\enddemo

\proclaim{Corollary}
Let\/ $A$ and\/ $B$ be two finite sets$.$
Let\/ $f$ be a function from\/ $A$ to\/ $\eusm P(B)$ and
let\/ $w$ be a function from\/ $B$ to\/ $\Bbb N\cup\{0\}.$
Let\/ a function\/ $F\:\eusm P(A)\to\eusm P(B)$ be defined by\/
$F(X)=\bigcup_{x\in X}f(x).$
Then the following are equivalent\/$:$
\medskip

\itemitem{$(a)$}
There exists a function\/ $h\:A\to B$ such that\/$:$
\itemitemitem{$(1)$}
for all\/ $x\in A,$ $h(x)\in f(x),$ and
\itemitemitem{$(2)$}
for each\/ $y\in B,$ the full pre-image of\/ $y$ under\/ $h$
consists of at most\/ $w(y)$ elements$.$
\smallskip

\itemitem{$(b)$}
For each subset\/ $X$ of\/ $A,$ $\displaystyle \|X\|\le\sum_{y\in F(X)}w(y).$

\endproclaim

\demo{Proof\rm}
The implication $(a)\Rightarrow(b)$ is obvious.
To prove the converse implication, consider the set $B'$
and functions $f'\:A\to\eusm P({B'})$ and $F'\:\eusm P(A)\to\eusm P({B'})$
defined by the formulae
$$
\aligned
B'&=
\bigl\{\,{(} b,n{)}\bigm| b\in B,\ n\in\Bbb N\cup\{0\},\ n<w(b)\,\bigr\},\\
f'(x)&=
\bigl\{\,{(} b,n{)}\bigm| b\in f(x),\ n\in\Bbb N\cup\{0\},\ n<w(b)\,\bigr\},\\
F'(X)&=\bigcup\{\,f'(x)\mid x\in X\,\}.
\endaligned
$$
Since
$$
\Bigl(\forall X\subset A\Bigr)
\Bigl(\|F'(X)\|=\sum_{y\in F(X)}w(y)\Bigr),
$$
it follows from the previous lemma that there exists an injection
$h'\:A\to B'$ such that
$$
\bigl(\forall x\in A\bigr)\bigl(h'(x)\in f'(x)\bigr).
$$
Let $h=p_1\kl\circ h'$ where $p_1\kl$ is the first projection from
$B\times\bigl(\Bbb N\cup\{0\}\bigr)$ onto~$B$.
The function $h$ is a desired one.
\qed
\enddemo

\proclaim{Estimating Lemma 1}
Let\/ $\Delta$ be an S-map with\/ $c\le3$ contours
satisfying the condition\/ $\Cal Z(2).$
Let\/ $k_1\kl,$ $k_2\kl,$ $k_3\kl$ be functions\/ $\Delta_\kl(2)\to\Bbb N\cup\{0\}.$
Suppose\/ $\Delta$ satisfies the conditions\/
$\Cal A_1\kl(k_1\kl),$ $\Cal A_2\kl(k_2\kl),$ $\Cal A_3\kl(k_3\kl).$
Suppose\/ $\Delta$ has at least one double-selected internal edge$.$
Let\/ $U$ be a non-overlapping set of double-selected internal arcs of\/ $\Delta$
covering the set of all double-selected internal edges of\/ $\Delta,$
with the minimal possible number of elements$.$
Then
$$
\|U\|\le\sum_{f\in\Delta(2)}\bigl(3+k_1\kl(f)+k_2\kl(f)+k_3\kl(f)\bigr)+2c-6,
$$
and there exists a function\/ $h\:U\to\Delta(2)$ such that\/$:$
\medskip

\itemitem{$(1)$}
each arc\/ $u\in U$ is incident to the face\/ $h(u),$ and
\smallskip

\itemitem{$(2)$}
the full pre-image of each face\/ $f\in\Delta(2)$ under\/ $h$
consists of at most\/ $3+k_1\kl(f)+k_2\kl(f)+k_3\kl(f)$ elements$.$

\endproclaim

\demo{Proof\rm}
For each $x\in U$, let $g(x)$ be the set of all faces of $\Delta$ incident to~$u$.
For each $X\subset U$, let $G(X)=\bigcup_{x\in X}g(x)$.
Take an arbitrary nonempty subset $E$ of~$U$.
It is to be proved that
$$
\|E\|\le\sum_{f\in G(E)}\bigl(3+k_1\kl(f)+k_2\kl(f)+k_3\kl(f)\bigr)+2c-6.
$$

Let $\bar\Delta$ be a spherical closure of~$\Delta$.
Without loss of generality, assume that $E$ cannot be split into
the disjoint union of two nonempty subsets $E_1\kl$ and $E_2\kl$ so that
$G(E_1\kl)\cap G(E_2\kl)=\emptyset$.
Then there exists an estimating $2$-complex $\Phi$ for $E$ in~$\bar\Delta$.

Denote $G(E)$ by~$F$.
Let $V$ be the set of all connected components of the $2$-complex
obtained from $\bar\Delta$ by removing all the faces that belong to $F$,
and all the edges and intermediate vertices of all the arcs
that belong to~$E$.
Let $W$ denote the set of all components of this complex that contain
improper faces of~$\bar\Delta$.
Note that $\|W\|\le c$.
Let $\hat E$ be the set of all oriented arcs associated with
arcs from~$E$.
Let
$$
\alpha_0\Kl\:V\to\Phi(0),\quad
\alpha_1'\:\hat E\to\hat\Phi(1),\quad
\alpha_2\Kl\:F\to\Phi(2)
$$
be such bijections as in the definition of an estimating $2$-complex.
Let $W'$ be the $\alpha_0\kl$-image of the set~$W$.

Consider $\Phi$.
Let $n_0\kl=\|\Phi(0)\|=\|V\|$, $n_1\kl=\|\Phi(1)\|=\|E\|$, $n_2\kl=\|\Phi(2)\|=\|F\|$.
Let $l=\|W'\|=\|W\|$.
Consider the vertices of $\Phi$ that are not in~$W'$.
Let $m_0^{(1)}$ be the number of such vertices of degree $1$,
$m_0^{(2)}$ be the number of such vertices of degree $2$,
$m_0^{(\ge3)}$ be the number of such vertices of degree at least~$3$.

On one hand, by Euler's Formula,
$$
n_0\kl-n_1\kl+n_2\kl=2;
$$
on the other hand,
$$
2n_1\Kl=\sum_{x\in\Phi(0)}d_{\Phi}\Kl(x)\ge m_0^{(1)}+2m_0^{(2)}+3m_0^{(\ge3)}+\sum_{w'\in W'}d_{\Phi}\Kl(w').
$$
Therefore,
$$
\aligned
3n_0\Kl-3n_1\Kl+3n_2\Kl+2n_1\Kl
&\ge6+m_0^{(1)}+2m_0^{(2)}+3m_0^{(\ge3)}+l,\\
3m_0^{(1)}+3m_0^{(2)}+3m_0^{(\ge3)}+3l-n_1\Kl+3n_2\Kl
&\ge6+m_0^{(1)}+2m_0^{(2)}+3m_0^{(\ge3)}+l,\\
3n_2\Kl+2m_0^{(1)}+m_0^{(2)}+2l-6&\ge n_1\Kl.
\endaligned
$$
In  other terms,
$$
\|E\|\le3\|F\|+2m_0^{(1)}+m_0^{(2)}+2l-6.
$$
Note that $2l-6\le2c-6\le0$.
Now, estimate the summand $2m_0^{(1)}+m_0^{(2)}$.

Let $M(f)$ be the set of all maximal selected subpaths
of the contour of a face $f$ of~$\Delta$.
Let $M_2\kl(f)$ be the set of all such $p_1\kl\in M(f)$ that
there exist $p_2\kl\in M(f)$ and a nontrivial path $q$ such that
$p_1\kl qp_2\kl$ is a subpath of the contour of $f$,
and no subpath of $q$ is selected.
Let $M_3\kl(f)$ be the set of all such $p_1\kl\in M(f)$ that
there exists $p_2\kl\in M(f)$ such that
$p_1\kl p_2\kl$ is a non-reduced subpath of the contour of~$f$.
The conditions $\Cal A_1\kl(k_1\kl)$, $\Cal A_2\kl(k_2\kl)$, $\Cal A_3\kl(k_3\kl)$
imply that for every $f\in\Delta(2)$,
$$
\|M(f)\|\le k_1\kl(f),\quad
\|M_2\kl(f)\|\le k_2\kl(f),\quad
\|M_3\kl(f)\|\le k_3\kl(f).
$$

Let
$$
N=\bigcup_{f\in F}M(f),\quad
N_2\kl=\bigcup_{f\in F}M_2\kl(f),\quad
N_3\kl=\bigcup_{f\in F}M_3\kl(f).
$$
Define a function $s\:N\to V$ by the following rule:
$s(p)$ is the last element of $V$ that the path $p$ meets,
\ie, such an element of $V$ that some terminal subpath of $p$
has a vertex in $s(p)$ and has no vertices
in any other element of~$V$.
It is easy to see that $s$ is well-defined.
Assign a weight to every element of $N$:
let the elements of $N_2\kl\sqcup N_3\kl$ have weight $2$,
and the other elements of $N$ have weight~$1$.
Let the weight of every subset $X$ of $N$ be the sum of the
weights of all its elements.
Let the weight of every $v\in V$ be the weight of its full
pre-image with respect to~$s$.
Let the weight of every $v'\in\Phi(0)$ be the weight of~$\alpha_0^{-1}(v')$.
It is to be proved now that every vertex of $\Phi$ of degree $2$
which is not in $W'$ weighs at least $1$,
and every vertex of $\Phi$ of degree $1$
which is not in $W'$ weighs at least~$2$.

Let $v'\in\Phi(0)\setminus W'$ and $v=\alpha_0^{-1}(v')$.
View the $2$-complex $v\in V$ as a submap of~$\bar\Delta$.
It is a disc map because of the assumptions about~$E$.

Consider the case $d_{\Phi}\kl(v')=2$.
It is enough to show that the full pre-image of $v$ under $s$ is not empty.
Suppose it is empty.
Let $u_1\kl$, $u_2\kl$ be those oriented arcs from $\hat E$ whose
terminal vertices are in $v$ (there are exactly $2$ such oriented arcs).
Let $f_i\kl$ be the face (from $F$) whose contour has $u_i\kl$ as a subpath, $i=1$,~$2$.
Then there exist subpaths $q_1\kl$, $q_2\kl$ of the contours
of $f_1\kl$, $f_2\kl$ respectively such that for $i=1$,~$2$,
the path $u_i\Kl q_i\Kl u_{3-i}^{-1}$ is a subpath of the contour of~$f_i\kl$.
Let $q_1\kl$, $q_2\kl$ be the shortest such paths.
Note that the cyclic path $q_1\kl q_2\kl$ is a representative
of the contour of the map~$v$.
Since the full pre-image of $v$ under $s$ is empty,
both paths $u_1\Kl q_1\Kl u_{2}^{-1}$ and $u_2\Kl q_2\Kl u_{1}^{-1}$
are selected.
Indeed, suppose that for example the path $u_1\Kl q_1\Kl u_{2}^{-1}$ is not selected.
Then there exists a maximal selected subpath $p$ of the contour of $f_1\kl$
that contains $u_1\kl$ as a subpath.
The path $p$ is an element of $M(f_1\kl)$.
Since $p$ cannot contain $u_1\Kl q_1\Kl u_{2}^{-1}$ as a subpath,
it follows that $s(p)=v$, which contradicts to the above assumption.
Hence, the paths $u_1\Kl q_1\Kl u_{2}^{-1}$ and $u_2\Kl q_2\Kl u_{1}^{-1}$
are selected; therefore, they are reduced.

Suppose the map $v$ is degenerate (for example, consists just of one vertex).
Then the paths $u_1\Kl q_1\Kl u_{2}^{-1}$ and $u_2\Kl q_2\Kl u_{1}^{-1}$
are mutually inverse.
They cannot be oriented arcs because it would contradict to the minimality
of the number of elements of $U$---two arcs corresponding to $u_1\kl$ and $u_2\kl$
could be replaced with the one corresponding to $u_1\Kl q_1\Kl u_{2}^{-1}$.
Therefore, $u_1\Kl=u_{2}^{-1}$, $v$ is the only element of $V$, and
$\Delta$ is exceptional.
Then $N$ is the full pre-image of $v$ under $s$, which is empty.
Therefore, all subpaths of the contours of $f_1\kl$ and $f_2\kl$ are selected.
Hence, $\Delta$ does not satisfy the condition $\Cal Z(2)$---not even $\Cal Z(0)$.
This gives a contradiction.

Suppose the map $v$ has at least one face.
Let $S$ be the set of all nontrivial paths in the set $\{q_1\kl,q_2\kl\}$.
The set $S$ consists of selected paths and covers the contour of the map~$v$.
Hence, by the Corollary of Proposition~3.1,
$\Delta$ cannot satisfy the condition $\Cal Z(2)$.
Contradiction.
Thus, if $d_{\Phi}\kl(v')=2$, the weight of $v'$ is at least~$1$.

Consider the case $d_{\Phi}\kl(v')=1$.
Let $u$ be the oriented arcs from $\hat E$ whose
terminal vertex is in $v$ (there is exactly $1$ such oriented arc).
Let $f$ be the face (from $F$) whose contour has $u$ as a subpath.
Then there exists a path $q$ such that
$uqu^{-1}$ is a subpath of the contour of~$f$.
Let $q$ be the shortest such path.
Note that the cyclic path $q$ is a representative
of the contour of the map~$v$.

Suppose the map $v$ is degenerate.
Then the path $uqu^{-1}$ is not reduced.
Therefore, there exists at least
one element of $M_2\kl(f)$ or $M_3\kl(f)$ in the $s$-pre-image of $v$,
and the weight of $v$ is at least~$2$.

Suppose the map $v$ has at least one face.
Then $q$ is nontrivial.
Since $\Delta$ satisfies the condition $\Cal Z(2)$,
it follows from the Corollary of Proposition~3.1 that
the path $q$ is not selected.
Moreover, it cannot be the product of two
selected subpaths.
Therefore, the full pre-image of $v$ under $s$
contains either an element of $M_2\kl(f)\cup M_3\kl(f)$
or at least two elements of $M(f)$.
In either case, the weight of $v$ is at least~$2$.

On one hand, the sum of the weights of all vertices of $\Phi$
is at least $2m_0^{(1)}+m_0^{(2)}$;
on the other hand, it is not greater than $\|N\|+\|N_2\kl\|+\|N_3\kl\|$.
Therefore,
$$
2m_0^{(1)}+m_0^{(2)}\le\sum_{f\in F}\bigl(k_1\Kl(f)+k_2\Kl(f)+k_3\Kl(f)\bigr).
$$

Thus, for an arbitrary nonempty subset $E$ of $U$,
$$
\aligned
\|E\|&\le3\|G(E)\|+\sum_{f\in G(E)}\bigl(k_1\kl(f)+k_2\kl(f)+k_3\kl(f)\bigr)+2c-6\\
&\qquad=\sum_{f\in G(E)}\bigl(3+k_1\kl(f)+k_2\kl(f)+k_3\kl(f)\bigr)+2c-6.
\endaligned
$$
In particular,
$$
\|U\|\le\sum_{f\in\Delta(2)}\bigl(3+k_1\kl(f)+k_2\kl(f)+k_3\kl(f)\bigr)+2c-6.
$$

By the Corollary of Hall's lemma,
there exists a function $h\:U\to\Delta(2)$ such that
each $u\in U$ is incident to the face $h(u)$, and
the full pre-image of each face $f\in\Delta(2)$ under $h$
consists of at most $3+k_1\kl(f)+k_2\kl(f)+k_3\kl(f)$ elements.
\qed
\enddemo

\proclaim{Estimating Lemma 2}
Let\/ $\Delta$ be a non-degenerate map$.$
Let\/ $S$ be the set of all\/ $($non-ordered\/$)$ pairs of distinct
contiguous faces of\/~$\Delta.$
Then\/ $\|S\|<3\|\Delta(2)\|,$
and there exists a function\/ $h\:S\to\Delta(2)$ such that\/$:$
\medskip

\itemitem{$(1)$}
each pair\/ $P\in S$ contains the face\/ $h(P),$ and
\smallskip

\itemitem{$(2)$}
the full pre-image of each face\/ $f\in\Delta(2)$ under\/ $h$
consists of at most\/ $3$ pairs$.$

\endproclaim

\demo{Proof\rm}
Define functions $g\:S\to\eusm P({\Delta(2)})$ and
$G\:\eusm P(S)\to\eusm P({\Delta(2)})$ by the formulae:
$$
g(P)=P,\qquad
G(X)=\bigcup X.
$$
Take an arbitrary nonempty subset $X$ of~$S$.
Let $F=\bigcup X=G(X)$.
The set $F$ is nonempty.
Let $\Gamma$ be the contiguity graph for $F$ in~$\Delta$.
Graph $\Gamma$ is planar, has no loops and no multiple edges.
Obviously, $\|\Gamma(0)\|=\|F\|$ and $\|\Gamma(1)\|\ge\|X\|$.
By Proposition~2.1, $\|X\|<3\|F\|$.

Thus, for any nonempty subset $X$ of $S$, $\|X\|<3\|G(X)\|$.
In particular,
$$
\|S\|<3\|\Delta(2)\|.
$$
Now, the conclusion easily follows from the Corollary of Hall's lemma.
\qed
\enddemo

\head
6. Statement of The Main Theorem
\endhead

\definition{Definition}
Let $\Delta$ be a simple S-map.
Let $n_1\kl$ be the number of edges of $\Delta$,
$n_1^{(e)}$ be the number of external edges of $\Delta$, and
$S$ be the number of selected external edges of~$\Delta$.
Let $\mu$ be a real number.
The S-map $\Delta$ is said to satisfy the {\it condition\/}
$\Cal X(\mu)$ if
$$
S\ge n_1\Kl-\mu(2n_1\Kl-n_1^{(e)})=\|\Delta(1)\|-\mu\sum_{f\in\Delta(2)}d(f).
$$
\enddefinition

\remark{Remark\/ \rm6.1}
In the above notation, if $\mu<1$ and $\Delta$ satisfies
the condition $\Cal X(\mu)$, then
$$
S\ge n_1\Kl-\mu\bigl(2n_1\Kl-n_1^{(e)}\bigr)
\ge\max\left\{\Bigl(1-\frac{\mu}{1-\mu}\Bigr)n_1\Kl,\,
\bigl(1-2\mu\bigr)\bigl(2n_1\Kl-n_1^{(e)}\bigr)\right\}.
$$
Indeed, since $n_1^{(e)}\ge S$ and $\mu<1$,
$$
\aligned
n_1^{(e)}&\ge n_1\Kl-\mu\bigl(2n_1\Kl-n_1^{(e)}\bigr),\\
n_1\Kl&\ge2n_1\Kl-n_1^{(e)}-\mu\bigl(2n_1\Kl-n_1^{(e)}\bigr),\\
n_1\Kl&\ge(1-\mu)\bigl(2n_1\Kl-n_1^{(e)}\bigr),\\
\frac{n_1\Kl}{1-\mu}&\ge2n_1\Kl-n_1^{(e)}.
\endaligned
$$
Therefore,
$$
n_1\Kl-\mu\bigl(2n_1\Kl-n_1^{(e)}\bigr)\ge(1-2\mu)\bigl(2n_1\Kl-n_1^{(e)}\bigr)
$$
and
$$
n_1\Kl-\mu\bigl(2n_1\Kl-n_1^{(e)}\bigr)\ge\Bigl(1-\frac{\mu}{1-\mu}\Bigr)n_1\Kl.
$$
\endremark

\remark{Remark\/ \rm6.2}
If a simple S-map $\Delta$ satisfies the condition $\Cal X(\mu)$
for some $\mu<\frac{1}{2}$, then $\Delta$ is not spherical.
\endremark

\proclaim{Main Theorem}
Let\/ $\Delta$ be a simple S-map with at most\/ $3$ contours$.$
Let
$$
\gathered
k_1\kl,k_2\kl,k_3\kl\:\Delta(2)\to\Bbb N\cup\{0\},\quad
\lambda_1\kl,\lambda_2\kl,\lambda_3\kl,\lambda_4\kl\:\Delta(2)\to[0,1],\\
\bar a={(} k_1\kl,k_2\kl,k_3\kl;\lambda_1\kl,\lambda_2\kl,\lambda_3\kl,\lambda_4\kl{)}.
\endgathered
$$
Let\/ $\mu$ be a real number$.$
Suppose\/ $\Delta$ satisfies the condition\/ $\Cal A(\bar a).$
Suppose
$$
\gathered
2\cdot\max_{\Delta(2)}\bigl(\lambda_1\kl+(3+k_1\kl+k_2\kl+k_3\kl)\lambda_2\kl+3\lambda_3\kl\bigr)
+\max_{\Delta(2)}\bigl((2+k_1\kl)\lambda_2\kl+2\lambda_4\kl\bigr)<1,\\
\max_{\Delta(2)}\bigl(\lambda_1\kl+(3+k_1\kl+k_2\kl+k_3\kl)\lambda_2\kl+3\lambda_3\kl\bigr)\le\mu.
\endgathered
$$
Then the S-map\/ $\Delta$ satisfies the conditions\/ $\Cal X(\mu)$ and\/ $\Cal Z(2).$
\endproclaim

The main idea of the proof is to verify the conditions $\Cal X(\mu)$ and $\Cal Z(2)$
by simultaneous induction on the number of internal edges of~$\Delta$.
Necessary inductive lemmas are proved in the next two sections.

\remark{Remark\/ \rm6.3}
Provided this theorem holds,
the condition $\Cal A(\bar a)$ and the inequalities
in the hypotheses of the theorem may be replaced with
the conjunction of only five conditions
$\Cal A_1\kl(k_1\kl)$, $\Cal A_4\kl(\lambda_1\kl)$, $\Cal A_5\kl(\lambda_2\kl)$,
$\Cal A_6\kl(\lambda_2\kl,\lambda_3\kl)$, $\Cal A_7\kl(\lambda_2\kl,\lambda_4\kl)$,
and with the inequalities
$$
\gathered
2\cdot\max_{\Delta(2)}\bigl(\lambda_1\kl+(3+2k_1\kl)\lambda_2\kl+3\lambda_3\kl\bigr)
+\max_{\Delta(2)}\bigl((2+k_1\kl)\lambda_2\kl+2\lambda_4\kl\bigr)<1,\\
\max_{\Delta(2)}\bigl(\lambda_1\kl+(3+2k_1\kl)\lambda_2\kl+3\lambda_3\kl\bigr)\le\mu,
\endgathered
$$
respectively.
Indeed, for any S-map $\Delta$ and function $k_1\kl\:\Delta(2)\to\Bbb N\cup\{0\}$
such that $\Delta$ satisfies $\Cal A_1\kl(k_1\kl)$,
there obviously exist $k_2\kl,k_3\kl\:\Delta(2)\to\Bbb N\cup\{0\}$ such that
$k_2\kl+k_3\kl\le k_1\kl$ point-wise, and $\Delta$ satisfies $\Cal A_2\kl(k_2\kl)$, $\Cal A_3\kl(k_3\kl)$.
Thus, $k_1\kl$ may be substituted for $k_2\kl+k_3\kl$ in the former inequalities
since $k_2\kl$ and $k_3\kl$ do not appear anywhere else in the statement of the theorem.
\endremark

\remark{Remark\/ \rm6.4}
In the case when all the functions
$k_1\kl$, $k_2\kl$, $k_3\kl$, $\lambda_1\kl$, $\lambda_2\kl$, $\lambda_3\kl$, $\lambda_4\kl$
are constant, the first inequality in the hypotheses is equivalent to
$$
2\lambda_1\kl+(8+3k_1\kl+2k_2\kl+2k_3\kl)\lambda_2\kl+6\lambda_3\kl+2\lambda_4\kl<1,
$$
and can be replaced with
$$
2\lambda_1\kl+(8+5k_1\kl)\lambda_2\kl+6\lambda_3\kl+2\lambda_4\kl<1.
$$
\endremark

\head
7. Inductive Lemma 1
\endhead

\proclaim{Inductive Lemma 1}
Let\/ $\Delta$ be an S-map$.$
Let
$$
k_1\kl\:\Delta(2)\to\Bbb N\cup\{0\},\quad
\lambda_2\kl,\lambda_4\kl\:\Delta(2)\to[0,1].
$$
Suppose\/ $\Delta$ satisfies the conditions\/
$\Cal A_1\kl(k_1\kl)$ and\/ $\Cal A_7\kl(\lambda_2\kl,\lambda_4\kl).$
Let\/ $\mu$ be a real number such that
the following inequality holds point-wise\/ $($face-wise\/$)$$:$
$$
2\mu+(2+k_1\kl)\lambda_2\kl+2\lambda_4\kl<1.
$$
Suppose that every proper simple disc S-submap of\/ $\Delta$
satisfies the condition\/ $\Cal X(\mu).$
Then\/ $\Delta$ satisfies the condition\/ $\Cal Z(2).$
\endproclaim

\demo{Proof\rm}
Suppose that such $\Delta$ does not satisfy the condition $\Cal Z(2)$.
Then in $\Delta$ there exists a simple disc submap $\Delta_0\kl$
such that there exists a set of selected paths with at most $2$
maximal elements which covers the contour of~$\Delta_0\kl$.
Let $\Delta_0\kl$ be such a submap with the minimal
possible number of internal edges.
View $\Delta_0\kl$ as an S-submap.
Let $c$ be the contour of~$\Delta_0\kl$.
Let $Q$ be a set consisting either of $1$ selected cyclic path $q_1\kl$ which is
a representative of $c$, or of $2$ selected paths $q_1\kl$, $q_2\kl$
whose product $q_1\kl q_2\kl$ is a representative of $c$ (clearly, such $Q$ exists).
Let $s=\|Q\|$.

The map $\Delta_0\kl$ is a proper submap of~$\Delta$.
Indeed, the edges that lie on the contour of $\Delta_0\kl$
are external in $\Delta_0\kl$ but internal in~$\Delta$.
Hence, by the assumptions,
the S-map $\Delta_0\kl$ satisfies the condition $\Cal X(\mu)$.

Let $f_0\kl$ be a face of $\Delta_0\kl$ which has at least
$(1-2\mu)d(f_0\kl)$ selected external edges in~$\Delta_0\kl$.
Such a face exists, otherwise
the number $S$ of selected external edges of $\Delta_0\kl$
would satisfy the inequality
$$
S<\sum_{f\in\Delta_0\kl(2)}(1-2\mu)d(f),
$$
in contradiction with the condition $\Cal X(\mu)$.

Let $P$ be a non-overlapping set of arcs of $\Delta$ such that
every element of $P$ is a selected external arc of $\Delta_0\kl$
(\ie, external and selected in $\Delta_0\kl$)
incident to $f_0\kl$ and lying on a path from $Q$,
and such that $P$ covers the set of all selected external edges of $\Delta_0\kl$
incident to $f_0\kl$;
moreover, let $P$ be such a set with the minimal possible number of elements.
For $i=1$,~$s$, let $P_i\kl$ be the set of those elements of $P$
that lie on~$q_i\kl$.
Let $m=\|P\|$, let $m_i\kl=\|P_i\kl\|$, $i=1$,~$s$.
For $i=1$,~$s$, enumerate all the subpaths of $q_i^{-1}$
which are oriented arcs corresponding to elements of $P_i\kl$,
according to the order in which they appear on $q_i^{-1}$:
$p_1^{(i)}$, \dots, $p_{m_i\kl}^{(i)}$.
Check now that for arbitrary $i$, $j$,
there is no selected subpath of the contour of $f_0\kl$
with initial subpath $p_j^{(i)}$ and terminal subpath~$p_{j+1}^{(i)}$.

Suppose on the contrary that for some $i\in\{1,s\}$ and
some $j\in\{\,n\in\Bbb N\mid 1\le n\le m_i\kl-1\,\}$,
there is a selected subpath of the contour of $f_0\kl$ with initial subpath
$p_j^{(i)}$ and terminal subpath~$p_{j+1}^{(i)}$.
Let $p_j^{(i)}xp_{j+1}^{(i)}$ be the shortest such selected subpath.
Then the path $x$ is either selected or trivial;
in either case, it is reduced.
Let $y$ be the subpath of $q_i\kl$ that
starts at the initial vertex of~$p_{j+1}^{(i)}$ and ends
at the terminal vertex of~$p_j^{(i)}$.
The path $y$ is reduced as a subpath of the contour of~$\Delta_0\kl$.
Clearly, the cyclic path $xy$ is a representative of the
contour of some disc submap $\Phi$ of the map~$\Delta_0\kl$.
Since the path $p_{j+1}^{(i)-1}yp_j^{(i)-1}$ is a subpath of $q_i\kl$,
and the path $p_j^{(i)}xp_{j+1}^{(i)}$ is selected,
the path $x$ cannot be inverse to the path $y$,
otherwise it would contradict
to the minimality of the number of elements of~$P$.
Therefore, the map $\Phi$ is non-degenerate.

By Proposition~3.1, there exists a maximal simple disc submap of $\Phi$
whose contour may be covered by $1$ or $2$ paths each of which
is a nontrivial subpath of either $x$ or~$y$.
Let $\Phi_1\kl$ be such a submap of~$\Phi$.
The map $\Phi_1\kl$ is also a submap of~$\Delta$.
The contour of $\Phi_1\kl$ is covered by a set of $1$ or $2$ selected paths.
The map $\Phi_1\kl$ clearly has fewer internal edges than~$\Delta_0\kl$.
This contradicts to the choice of
$\Delta_0\kl$---to the minimality of its number of internal edges.

Since for any $i$ and $j$, no subpath of the contour of $f_0\kl$
with initial subpath $p_j^{(i)}$ and terminal subpath $p_{j+1}^{(i)}$
is selected,
the number of maximal selected subpaths of the contour of $f_0\kl$
is not less than $m_1\kl-1=m-1$ in the case $s=1$,
and not less than $m_1\kl+m_2\kl-2=m-2$ in the case $s=2$.
Hence, in both cases, $m\le k_1\kl(f_0\kl)+2$
by the condition $\Cal A_1\kl(k_1\kl)$.

Let
$$
S_i\kl=\sum_{u\in P_i\kl}|u|=\sum_{j=1}^{m_i\kl}|p_j^{(i)}|,\quad i=1, s,
\qquad
S=\sum_{u\in P}|u|=\sum_{i=1}^sS_i\kl.
$$
Note that $S$ equals the number of selected external edges
of $f_0\kl$ in~$\Delta_0\kl$.
For each $i=1$,~$s$, the condition
$\Cal A_7\kl(\lambda_2\kl,\lambda_4\kl)$ for $\Delta$ implies that
$$
S_i\kl\le\bigl(m_i\kl\lambda_2\kl(f_0\kl)+\lambda_4\kl(f_0\kl)\bigr)d(f_0\kl).
$$
Therefore, on one hand,
$$
S\le\Bigl(m\lambda_2\kl(f_0\kl)+s\lambda_4\kl(f_0\kl)\Bigr)d(f_0\kl)
\le\Bigl(\bigl(2+k_1\kl(f_0\kl)\bigr)\lambda_2\kl(f_0\kl)+2\lambda_4\kl(f_0\kl)\Bigr)d(f_0\kl),
$$
on the other hand,
$$
S\ge(1-2\mu)d(f_0\kl).
$$
Hence,
$$
\bigl(2+k_1\kl(f_0\kl)\bigr)\lambda_2\kl(f_0\kl)+2\lambda_4\kl(f_0\kl)\ge1-2\mu,
$$
and this contradicts to the inequality
$$
2\mu+(2+k_1\kl)\lambda_2\kl+2\lambda_4\kl<1.\hbox to 0pt{\qed\hss}
$$
\enddemo

\head
8. Inductive Lemma 2
\endhead

\proclaim{Inductive Lemma 2}
Let\/ $\Delta$ be a simple S-map with at most\/ $3$ contours$.$
Let
$$
k_1\kl,k_2\kl,k_3\kl\:\Delta(2)\to\Bbb N\cup\{0\},\quad
\lambda_1\kl,\lambda_2\kl,\lambda_3\kl\:\Delta(2)\to[0,1].
$$
Suppose\/ $\Delta$ satisfies the conditions\/
$\Cal Z(2),$ $\Cal A_1\kl(k_1\kl),$ $\Cal A_2\kl(k_2\kl),$ $\Cal A_3\kl(k_3\kl),$
$\Cal A_4\kl(\lambda_1\kl),$ $\Cal A_5\kl(\lambda_2\kl),$ $\Cal A_6\kl(\lambda_2\kl,\lambda_3\kl).$
Let\/ $\mu$ be a real number such that
the following double inequality holds point-wise\/ $($face-wise\/$)$$:$
$$
\lambda_1\kl+(3+k_1\kl+k_2\kl+k_3\kl)\lambda_2\kl+3\lambda_3\kl\le\mu<1.
$$
Then\/ $\Delta$ satisfies the condition\/ $\Cal X(\mu).$
\endproclaim

\demo{Proof\rm}
Let $n_1\kl$ be the number of edges of $\Delta$,
$n_1^{(e)}$ be the number of external edges of $\Delta$,
$S$ be the number of selected external edges of~$\Delta$.
Note that since $\Delta$ is simple,
$$
\sum_{f\in\Delta(2)}d(f)=2n_1\Kl-n_1^{(e)}.
$$

First, estimate the number of all edges that are neither selected external,
nor double-selected internal.
Denote this number by~$S_1'$.
Since $S_1'$ is less than or equal to the number of non-selected oriented
edges of the contours of all the faces of $\Delta$,
and $\Delta$ satisfies the condition $\Cal A_4\kl(\lambda_1\kl)$,
$$
S_1'\le\sum_{f\in\Delta(2)}\lambda_1\Kl(f)d(f).
$$

Second, estimate the number of double-selected internal edges.
Denote this number by~$S_2'$.

Let $U$ be a minimal by the number of elements non-overlapping set
of double-selected internal arcs of $\Delta$
covering the set of all double-selected edge of~$\Delta$.
Let for any $f_1\kl,f_2\kl\in\Delta(2)$, $B(\{f_1\kl,f_2\kl\})$ be the
set of all arcs between $f_1\kl$ and $f_2\kl$ that are elements of $U$
(if $f_1\kl=f_2\kl$, then $B(\{f_1\kl,f_2\kl\})=B(\{f_1\kl\})$ is the set of all
elements of $U$ which are intro-facial arcs of~$f_1\kl$).
Let $M$ be the set of all such pairs or singletons $X$ of faces of $\Delta$
that $B(X)\ne\emptyset$.
Then
$$
S_2'=\sum_{u\in U}|u|=\sum_{X\in M}\sum_{u\in B(X)}|u|.
$$
Let $N$ be the set of all pairs of distinct contiguous faces of $\Delta$
that have at least one double-selected edge between them.

Let $h_1\kl\:U\to\Delta(2)$ and $h_2\kl\:N\to\Delta(2)$ be such functions that
\medskip

\itemitem{$(1)$}
each arc $u\in U$ is incident to the face~$h_1\kl(u)$;
\smallskip

\itemitem{$(2)$}
the full pre-image of each face $f$ of $\Delta$ under $h_1\kl$
consists of at most $3+k_1\kl(f)+k_2\kl(f)+k_3\kl(f)$ arcs;
\smallskip

\itemitem{$(3)$}
each pair $P\in N$ contains the face~$h_2\kl(P)$;
\smallskip

\itemitem{$(4)$}
the full pre-image of each face $f$ of $\Delta$ under $h_2\kl$
consists of at most $3$ pairs.
\medskip

\noindent
Such a function $h_1\kl$ exists by Estimating Lemma~1
if $U$ is not empty, and it obviously exists if $U$ is empty.
Such a function $h_2\kl$ exists by Estimating Lemma~2.
(The conditions $\Cal A_1\kl(k_1\kl)$, $\Cal A_2\kl(k_2\kl)$, $\Cal A_3\kl(k_3\kl)$,
and $\Cal Z(2)$ have been used here.)

Let $K$ be the subset of $\Delta(2)\times\Delta(2)$
consisting of all ordered pairs ${(} f_1\kl,f_2\kl{)}$
such that $\{f_1\kl,f_2\kl\}\in M$.
Define functions $r_1\kl\:K\to\Bbb R$ and $r_2\kl\:K\to\Bbb R$
as follows.
For any ${(} f_1\kl,f_2\kl{)}\in K$, let
$$
\aligned
r_1\Kl(f_1\Kl,f_2\Kl)&=\bigl\|h_1^{-1}(f_1\Kl)\cap B(\{f_1\Kl,f_2\Kl\})\bigr\|\lambda_2\Kl(f_1\Kl)d(f_1\Kl),\\
r_2\Kl(f_1\Kl,f_2\Kl)&=\bigl\|h_2^{-1}(f_1\Kl)\cap \{\{f_1\Kl,f_2\Kl\}\}\bigr\|\lambda_3\Kl(f_1\Kl)d(f_1\Kl).
\endaligned
$$
Let $r=r_1\kl+r_2\kl$.

Take arbitrary distinct faces $f_1\kl$ and $f_2\kl$ such that $\{f_1\kl,f_2\kl\}\in M$.
For the sake of definiteness, assume that $h_2\kl(\{f_1\kl,f_2\kl\})=f_1\kl$.
Then
$$
\aligned
&r(f_1\Kl,f_2\Kl)+r(f_2\Kl,f_1\Kl)\\
&\qquad=\bigl(r_1\Kl(f_1\Kl,f_2\Kl)+r_1\Kl(f_2\Kl,f_1\Kl)\bigr)
+\bigl(r_2\Kl(f_1\Kl,f_2\Kl)+r_2\Kl(f_2\Kl,f_1\Kl)\bigr)\\
&\qquad=\bigl\|h_1^{-1}(f_1\Kl)\cap B(\{f_1\Kl,f_2\Kl\})\bigr\|\lambda_2\Kl(f_1\Kl)d(f_1\Kl)
+\bigl\|h_1^{-1}(f_2\Kl)\cap B(\{f_1\Kl,f_2\Kl\})\bigr\|\lambda_2\Kl(f_2\Kl)d(f_2\Kl)\\
&\qquad\qquad+\lambda_3\Kl(f_1\Kl)d(f_1\Kl)\\
&\qquad\ge\bigl\|B(\{f_1\Kl,f_2\Kl\})\bigr\|\min\{\lambda_2\Kl(f_1\Kl)d(f_1\Kl),\lambda_2\Kl(f_2\Kl)d(f_2\Kl)\}
+\lambda_3\Kl(f_1\Kl)d(f_1\Kl).
\endaligned
$$
By the condition $\Cal A_6\kl(\lambda_2\kl,\lambda_3\kl)$,
$$
\aligned
\sum_{u\in B(\{f_1\kl,f_2\kl\})}|u|
&\le\bigl\|B(\{f_1\kl,f_2\kl\})\bigr\|\min\{\lambda_2\kl(f_1\kl)d(f_1\kl),\lambda_2\kl(f_2\kl)d(f_2\kl)\}
+\lambda_3\kl(f_1\kl)d(f_1\kl)\\
&\le r(f_1\kl,f_2\kl)+r(f_2\kl,f_1\kl).
\endaligned
$$

Now, take an arbitrary face $f$ such that $\{f\}\in M$.
Then
$$
\aligned
&r(f,f)=r_1\Kl(f,f)+r_2\Kl(f,f)\\
&\qquad=\bigl\|h_1^{-1}(f)\cap B(\{f\})\bigr\|\lambda_2\Kl(f)d(f)+0\\
&\qquad=\bigl\|B(\{f\})\bigr\|\lambda_2\Kl(f)d(f).
\endaligned
$$
By the condition $\Cal A_5\kl(\lambda_2\kl)$,
$$
\sum_{u\in B(\{f\})}|u|
\le\bigl\|B(\{f\})\bigr\|\lambda_2\kl(f)d(f)
=r(f,f).
$$

Thus,
$$
\aligned
S_2'&=\sum_{u\in U}|u|=\sum_{X\in M}\sum_{u\in B(X)}|u|\\
&\le\sum_{{(} f_1\kl,f_2\kl{)}\in K}r(f_1\Kl,f_2\Kl)
=\sum_{f_1\kl\in\Delta(2)}\,\sum_{f_2\kl:{(} f_1\kl,f_2\kl{)}\in K}r(f_1\Kl,f_2\Kl)\\
&=\sum_{f_1\kl\in\Delta(2)}
\biggl(\,\sum_{f_2\kl:{(} f_1\kl,f_2\kl{)}\in K}r_1\Kl(f_1\Kl,f_2\Kl)
+\sum_{f_2\kl:{(} f_1\kl,f_2\kl{)}\in K}r_2\Kl(f_1\Kl,f_2\Kl)\biggr)\\
&=\sum_{f_1\kl\in\Delta(2)}
\bigl(\|h_1^{-1}(f_1\Kl)\|\lambda_2\Kl(f_1\Kl)d(f_1\Kl)
+\|h_2^{-1}(f_1\Kl)\|\lambda_3\Kl(f_1\Kl)d(f_1\Kl)\bigr)\\
&\le\sum_{f_1\kl\in\Delta(2)}
\Bigl(\bigl(3+k_1\Kl(f_1\Kl)+k_2\Kl(f_1\Kl)+k_3\Kl(f_1\Kl)\bigr)\lambda_2\Kl(f_1\Kl)d(f_1\Kl)
+3\lambda_3\Kl(f_1\Kl)d(f_1\Kl)\Bigr)\\
&=\sum_{f\in\Delta(2)}
\Bigl(\bigl(3+k_1\Kl(f)+k_2\Kl(f)+k_3\Kl(f)\bigr)\lambda_2\Kl(f)
+3\lambda_3\Kl(f)\Bigr)d(f).
\endaligned
$$

Eventually,
$$
\aligned
S&=n_1\Kl-S_1'-S_2'\\
&\qquad\ge n_1\Kl-\sum_{f\in\Delta(2)}
\Bigl(\lambda_1\Kl(f)+\bigl(3+k_1\Kl(f)+k_2\Kl(f)+k_3\Kl(f)\bigr)\lambda_2\Kl(f)+3\lambda_3\Kl(f)\Bigr)d(f)\\
&\qquad\qquad\ge n_1\Kl-\sum_{f\in\Delta(2)}\mu d(f)=n_1\Kl-\mu\bigl(2n_1\Kl-n_1^{(e)}\bigr),
\endaligned
$$
since $\lambda_1\kl+(3+k_1\kl+k_2\kl+k_3\kl)\lambda_2\kl+3\lambda_3\kl\le\mu$.
\qed
\enddemo

\head
9. Proof of The Main Theorem
\endhead

\demo{Proof of the main theorem\rm}
Suppose the theorem is not true.
Then let $\Delta$ be a simple S-map,
let
$$
\gathered
k_1\kl,k_2\kl,k_3\kl\:\Delta(2)\to\Bbb N\cup\{0\},\quad
\lambda_1\kl,\lambda_2\kl,\lambda_3\kl,\lambda_4\kl\:\Delta(2)\to[0,1],\\
\bar a={(} k_1\kl,k_2\kl,k_3\kl;\lambda_1\kl,\lambda_2\kl,\lambda_3\kl,\lambda_4\kl{)},
\endgathered
$$
and let $\mu$ be a real number such that
they all together satisfy the hypotheses of the theorem and do not
satisfy the conclusion (\ie, provide a counterexample),
and such that the number of internal edges of $\Delta$
is the minimal possible under this condition
(\ie, the theorem holds whenever the S-map has fewer internal edges).
Without loss of generality, assume that
$$
\max_{\Delta(2)}\bigl(\lambda_1\kl+(3+k_1\kl+k_2\kl+k_3\kl)\lambda_2\kl+3\lambda_3\kl\bigr)=\mu.
$$
Then
$$
2\mu+\max_{\Delta(2)}\bigl((2+k_1\kl)\lambda_2\kl+2\lambda_4\kl\bigr)<1.
$$
In particular, $\mu<1$.

The hypotheses of Inductive Lemma~1 hold for $\Delta$, $\bar a$, and~$\mu$.
Therefore, by Inductive Lemma~1, $\Delta$ satisfies the condition $\Cal Z(2)$.
Therefore, by Inductive Lemma~2, $\Delta$ satisfies the condition $\Cal X(\mu)$.
Hence, the theorem holds for $\Delta$, $\bar a$, and $\mu$,
and this contradicts to how they were chosen.
\qed
\enddemo

\head
10. Lemma about exposed face
\endhead

\proclaim{Lemma about exposed face}
Let\/ $\Delta$ be a simple disc S-map$.$
Let
$$
k_1\kl\:\Delta(2)\to\Bbb N\cup\{0\},\quad
\lambda_2\kl,\lambda_4\kl\:\Delta(2)\to[0,1],\quad
\mu\in\Bbb R.
$$
Suppose\/ $\Delta$ satisfies the conditions\/
$\Cal Z(2),$ $\Cal A_1\kl(k_1\kl),$ $\Cal A_7\kl(\lambda_2\kl,\lambda_4\kl).$
Suppose every simple disc S-submap of\/ $\Delta$ satisfies
the condition\/ $\Cal X(\mu).$
Then there exists a face\/ $f$ in\/ $\Delta$ satisfying the following property\/$:$
if\/ $P$ is a non-overlapping set of selected external arcs incident to\/ $f,$
if $P$ covers the set of all selected external edges incident to\/ $f,$
and if the number of elements of\/ $P$ is the minimal possible under these assumptions$,$
then\/ $\|P\|\le k_1\kl(f)+1$ and
$$
\sum_{p\in P}|p|
\ge\Bigl(1-2\mu-\bigl(2+k_1\kl(f)\bigr)\lambda_2\kl(f)-\lambda_4\kl(f)+\|P\|\lambda_2\kl(f)\Bigr)d(f).
$$
$($Such a face\/ $f$ may be called\/ {\rm``}exposed\/{\rm''}$)$$.$
\endproclaim

\remark{Remark\/ \rm10.1}
It is possible to prove
a stronger statement than the one claimed in this lemma.
Namely, under the hypotheses of the lemma,
either there exists a face $f$ in $\Delta$
such that a corresponding set $P$ has at most $k_1\kl(f)+1$ elements,
and the total sum of their lengths is at least $(1-2\mu)d(f)$,
or there exist at least two distinct ``exposed'' faces
such as in the conclusion of the lemma.
\endremark

\remark{Remark\/ \rm10.2}
If $\Delta$, $\bar a={(} k_1\kl,k_2\kl,k_3\kl;\lambda_1\kl,\lambda_2\kl,\lambda_3\kl,\lambda_4\kl{)}$,
and $\mu$ satisfy the hypotheses of the main theorem,
then it follows from the main theorem that
$\Delta$, $k_1\kl$, $\lambda_2\kl$, $\lambda_4\kl$, $\mu$
satisfy the hypotheses of the lemma about exposed face.
\endremark

\demo{Proof\rm}
Case I:
there exists a simple disc submap $\Delta_0\kl$
of $\Delta$ whose contour is of the form $q_1\kl q_2\kl$
where $q_1\kl$ is a selected subpath of the contour of some face $f_0\kl$,
and $q_2\kl$ is a subpath of the contour of~$\Delta$.
Let $\Delta_0\kl$ be a minimal such map
(which do not contain any other such map as a submap).
View $\Delta_0\kl$ as an S-submap.
Let $q_1\kl$, $q_2\kl$ and $f_0\kl$ be as above.
It follows from the condition $\Cal Z(2)$ that $q_2\kl$ is nontrivial.

Pick a face $f_1\kl$ of $\Delta_0\kl$ which has at least $(1-2\mu)d(f_1\kl)$
selected external edges in~$\Delta_0\kl$.
Such a face exists because $\Delta_0\kl$ satisfies the condition $\Cal X(\mu)$.

Let $P_1\kl$ be a non-overlapping set of double-selected arcs between $f_1\kl$ and $f_0\kl$
such that $P_1\kl$ covers the set of all double-selected edges between $f_1\kl$ and~$f_0\kl$.
Let $P_2\kl$ be a non-overlapping set of selected external arcs incident to $f_1\kl$
such that $P_2\kl$ covers the set of all selected external edge incident to~$f_1\kl$.
Moreover, let $P_1\kl$ and $P_2\kl$ be such sets with the minimal possible numbers of elements.
Note that all the elements of $P_1\kl$ lie on the path $q_1\kl$
and all the element of $P_2\kl$ lie on the path~$q_2\kl$.
Let $m_1\kl=\|P_1\kl\|$, $m_2\kl=\|P_2\kl\|$.
Let $P=P_1\kl\sqcup P_2\kl$.
All the elements of $P$ are selected external arcs of $\Delta_0\kl$
incident to~$f_1\kl$.
By the choice of $f_1\kl$,
$$
\sum_{p\in P}|p|\ge(1-2\mu)d(f_1\kl).
$$
The goal is to prove that $m_2\kl\le k_1\kl(f_1\kl)+1$ and
$$
\sum_{p\in P_2\kl}|p|
\ge\Bigl(1-2\mu-\bigl(2+k_1\kl(f_1\kl)\bigr)\lambda_2\kl(f_1\kl)-\lambda_4\kl(f_1\kl)+m_2\kl\lambda_2\kl(f_1\kl)\Bigr)
d(f_1\kl);
$$
this will show that the face $f_1\kl$ is a desired one.

For $i=1$,~$2$, if $P_i\kl\ne\emptyset$, then
enumerate all the subpaths of $q_i^{-1}$
which are oriented arcs corresponding to elements of $P_i\kl$,
according to the order in which they appear on $q_i^{-1}$:
$p_1^{(i)}$, \dots, $p_{m_i\kl}^{(i)}$.
Check now that for arbitrary $i$, $j$,
there is no selected subpath of the contour of $f_1\kl$
with initial subpath $p_j^{(i)}$ and terminal subpath~$p_{j+1}^{(i)}$.

Suppose that for some $j\in\{\,n\in\Bbb N\mid 1\le n<m_1\kl\,\}$,
there is a selected subpath of the contour of $f_1\kl$ with initial subpath
$p_j^{(1)}$ and terminal subpath~$p_{j+1}^{(1)}$.
Let $p_j^{(1)}xp_{j+1}^{(1)}$ be the shortest such selected subpath.
Then the path $x$ is either selected or trivial;
in either case, it is reduced.
Let $y$ be the subpath of $q_1\kl$ that
starts at the initial vertex of~$p_{j+1}^{(1)}$ and ends
at the terminal vertex of~$p_j^{(1)}$.
The path $y$ is reduced as a subpath of the contour of~$\Delta_0\kl$.
Note that $p_{j+1}^{(1)-1}yp_j^{(1)-1}$ is a subpath of~$q_1\kl$.
The cyclic path $xy$ is a representative of the
contour of some disc submap $\Phi$ of~$\Delta$.
The path $x$ cannot be inverse to $y$, otherwise it would contradict
to the minimality of the number of elements of~$P_1\kl$.
Therefore, the map $\Phi$ is non-degenerate.
The Corollary of Proposition~3.1
gives a contradiction with the fact that $\Delta$ satisfies
the condition $\Cal Z(2)$.

Suppose that for some $j\in\{\,n\in\Bbb N\mid 1\le n<m_2\kl\,\}$,
there is a selected subpath of the contour of $f_1\kl$ with initial subpath
$p_j^{(2)}$ and terminal subpath~$p_{j+1}^{(2)}$.
Let $p_j^{(2)}xp_{j+1}^{(2)}$ be the shortest such selected subpath.
Since the disc map $\Delta$ is simple and the number of elements of $P_2\kl$
is minimal, the path $x$ cannot be trivial.
Therefore, the path $x$ is selected and reduced.
Let $y$ be the subpath of $q_2\kl$ that
starts at the initial vertex of~$p_{j+1}^{(2)}$ and ends
at the terminal vertex of~$p_j^{(2)}$.
The path $y$ is reduced as a subpath of the contour of~$\Delta_0\kl$.
Note that $p_{j+1}^{(2)-1}yp_j^{(2)-1}$ is a subpath of~$q_2\kl$.
The cyclic path $xy$ is a representative of the
contour of some disc submap $\Phi$ of~$\Delta$.
The path $x$ cannot be inverse to $y$, otherwise it would contradict
to the minimality of the number of elements of~$P_2\kl$.
Therefore, the map $\Phi$ is non-degenerate.
By Proposition~3.1, there exists a simple disc submap $\Delta_1\kl$
of $\Phi$ whose contour is covered by a set of $1$ or $2$ paths,
each of which is either a subpath of $x$ or a subpath of~$y$.
In any case, this leads to a contradiction either with the condition $\Cal Z(2)$,
or with the fact that $\Delta$ is simple, or with the choice of $\Delta_0$
(its minimality).

Thus, for arbitrary $i$ and $j$, there is no selected path
with initial subpath $p_j^{(i)}$ and terminal subpath~$p_{j+1}^{(i)}$.
Hence, it follows from the condition $\Cal A_1\kl(k_1\kl)$ that
$$
m_1\kl\le k_1\kl(f_1\kl)+1,\quad
m_2\kl\le k_1\kl(f_1\kl)+1,\quad
m_1\kl+m_2\kl\le k_1\kl(f_1\kl)+2.
$$

By the condition $\Cal A_7\kl(\lambda_2\kl,\lambda_4\kl)$,
$$
\aligned
\sum_{p\in P_1\kl}|p|&\le\Bigl(m_1\kl\lambda_2\kl(f_1\kl)+\lambda_4\kl(f_1\kl)\Bigr)d(f_1\kl)\\
&\qquad\le\Bigl(\bigl(k_1\kl(f_1\kl)+2-m_2\kl\bigr)\lambda_2\kl(f_1\kl)+\lambda_4\kl(f_1\kl)\Bigr)d(f_1\kl).
\endaligned
$$
Therefore,
$$
\aligned
\sum_{p\in P_2\kl}|p|
&\ge\Bigl(1-2\mu-\bigl(k_1\kl(f_1\kl)+2-m_2\kl\bigr)\lambda_2\kl(f_1\kl)-\lambda_4\kl(f_1\kl)\Bigr)d(f_1\kl)\\
&\qquad=\Bigl(1-2\mu-\bigl(2+k_1\kl(f_1\kl)\bigr)\lambda_2\kl(f_1\kl)-\lambda_4\kl(f_1\kl)
+m_2\kl\lambda_2\kl(f_1\kl)\Bigr)d(f_1\kl).
\endaligned
$$

Case II:
there exists no simple disc submap $\Delta_0\kl$
of $\Delta$ whose contour would be of the form $q_1\kl q_2\kl$
where $q_1\kl$ would be a selected subpath of the contour of some face,
and $q_2\kl$ would be a subpath of the contour of~$\Delta$.

Consider an arbitrary face $f$ of~$\Delta$.
Let $P$ be a non-overlapping set of selected external arcs incident to $f$
such that $P$ covers the set of all selected external edges incident to $f$,
and the number of elements of $P$ is the minimal possible.
Then $\|P\|\le k_1\kl(f)+1$.
Indeed, if $P$ is empty, this inequality is obvious.
Suppose $P$ is not empty.
For every selected subpath $q$ of the contour of $f$,
at most one element of $P$ lies on $q$
(because the number of elements of $P$ is minimal,
$\Delta$ is a simple disc map, $\Delta$ satisfies $\Cal Z(2)$,
and it is not Case~I).
If the set of all selected subpath of the contour of $f$
has no maximal element, then $P$ consists of one element
(because there exists a selected subpath $q$
of the contour of $f$ such that all the elements of $P$ lie on $q$),
and $\|P\|=1\le k_1\kl(f)+1$.
If the set of all selected subpath of the contour of $f$
has at least one maximal element, then every element of $P$
lies on some maximal selected subpath of the contour of $f$,
and $\|P\|\le k_1\kl(f)<k_1\kl(f)+1$ (by $\Cal A_1\kl(k_1\kl)$).

Now, pick a face $f_0\kl$ in $\Delta$ which has at least $(1-2\mu)d(f_0\kl)$
selected external edges.
The face $f_0\kl$ is a desired one.
\qed
\enddemo

\head
11. Group presentations and diagrams
\endhead

Recall that an (abstract) {\it group presentation\/}
is an ordered pair $\langle\,\goth A\,\|\,\Cal R\,\rangle$
where $\goth A$ is an arbitrary set, called {\it alphabet},
and $\Cal R$ is a set of words in the
{\it group alphabet\/}~$\goth A^{\pm1}$.
The elements of $\Cal R$ are called {\it defining words}.
The same group presentation may be also written as
$\langle\,\goth A\,\|\,R=1,\ R\in\Cal R\,\rangle$,
here the relations $R=1$, $R\in\Cal R$, are called {\it defining relations}.
It will be assumed in the rest of this paper
that $\Cal R$ does not contain the empty word.

A {\it group word\/} in the alphabet $\goth A$
is a word in the alphabet~$\goth A^{\pm1}$.
The {\it inverse\/} word to a group word $W$ is denoted by~$W^{-1}$.
The $m$th {\it power\/} of a group word $W$ for an integer $m$ is denoted by $W^m$
and defined as follows.
If $m$ is a positive integer, then $W^m$ is the result of concatenation
of $m$ copies of~$W$.
If $m$ is a negative integer, then $W^m$ is the result of concatenation
of $-m$ copies of~$W^{-1}$.
By definition, $W^0$ is the empty word.
If $A$ and $B$ are two group words, then define $A^B=BAB^{-1}$.
A group word is called {\it reduced\/} if
it has no subwords of the form $xx^{-1}$ where $x\in\goth A^{\pm1}$.
The reduced word obtained from a group word $W$ by cancelling one-by-one
all subwords of the form $xx^{-1}$,  $x\in\goth A^{\pm1}$,
is called the {\it reduced form\/} of $W$ (it is well-defined).
A group word is {\it cyclically reduced\/} if
it is reduced and its first letter is not inverse to its last one.
The empty word is cyclically reduced by definition.
Two word are {\it freely equal\/} if
their reduced forms coincide.
A set of words (or relations) is called {\it symmetrized\/} if
along with every word $W$ it contains all cyclic shifts of $W$ and~$W^{-1}$.
Sometimes, the set of defining words of a presentation will be assumed symmetrized.
A {\it cyclic word\/} is the set of all cyclic shifts of some word.
The cyclic word corresponding to a word $W$ (containing the word $W$)
will be denoted by~$\tilde W$.
The {\it length\/} of a cyclic word is the length of an arbitrary
representative of it.
A word $V$ is a {\it subword\/} of a cyclic word $\tilde W$ if
it is a subword of some representative of~$\tilde W$.

Let $\langle\,\goth A\,\|\,\Cal R\,\rangle$ be a group presentation.
Let $W$ be a group word in the alphabet~$\goth A$.
The relation $W=1$ is a {\it consequence\/} of the system of relations
$\{\,R=1\mid R\in\Cal R\,\}$ if%
, roughly speaking,
it may be deduced from relations of this system, associativity, and
cancellation of inverse letters.

Every presentation $\langle\,\goth A\,\|\,\Cal R\,\rangle$ defines
a group $G$, and there is a natural function from the set of all
group words in the alphabet $\goth A$ onto this group.
The image of a group word $W$ under this function is trivial in $G$ if and only if
the relation $W=1$ is a consequence of the system of relations
$\{\,R=1\mid R\in\Cal R\,\}$.
The images of the one-letter words associated with the elements of $\goth A$
generate the whole of~$G$.
Two group words $W_1\kl$ and $W_2\kl$ are said to be {\it equal in the group\/} $G$ if
their images in $G$ are equal, or, equivalently, if
the relation $W_1\kl=W_2\kl$ is a consequence of the relations
$R=1$, $R\in\Cal R$.

A {\it diagram\/} over a presentation $\langle\,\goth A\,\|\,\Cal R\,\rangle$
is a map together with a function which labels every oriented edge of this
map with a letter from $\goth A^{\pm1}$ so that mutually inverse oriented edges are labelled
with mutually inverse group letters, and for every face, the word that ``reads''
on some representative of its contour is an element of~$\Cal R^{\pm1}$.

In a diagram over $\langle\,\goth A\,\|\,\Cal R\,\rangle$,
let the {\it label\/} of an oriented edge be the element of $\goth A^{\pm1}$
that labels this oriented edge, the {\it label\/} of a path
be the group word that reads on this path,
the {\it label\/} of a cycle be the cyclic (group) word whose representatives
are the labels of representatives of this cycle,
the {\it label\/} of a face be the label of its contour.

According to van Kampen's Lemma (see \cite{O}),
a relation $W=1$ is a consequence of a system of relations $\{\,R=1\mid R\in\Cal R\,\}$
if and only if there exists a disc diagram $\Delta$ such that
the label of some representative of its contour is $W$,
and for every face $f$ of $\Delta$, the label of some representative
of the contour of $f$ either is an element of $\Cal R$, or is inverse to
some element of~$\Cal R$.
Such a diagram is called a {\it deduction diagram\/} for the word $W$
or relation $W=1$.

A pair of distinct faces $\{f_1\kl,f_2\kl\}$ in a diagram $\Delta$
is called {\it cancellable\/} if
there exist paths $p_1\kl$ and $p_2\kl$ with common initial oriented edge
and with equal labels such that $p_1\kl$ is a representative of the contour
of $f_1\kl$, and $p_2^{-1}$ is a representative of the contour of~$f_2\kl$.
A diagram without cancellable pairs of faces is called {\it reduced}.
If a relation $W=1$ is a consequence of a system of relations
$\{\,R=1\mid R\in\Cal R\,\}$, then, in fact, there exists a reduced deduction
diagram over $\langle\,\goth A\,\|\,\Cal R\,\rangle$ for the relation $W=1$.

\definition{Definition}
A symmetrized set $S$ of group words in a given alphabet satisfies the
{\it small cancellation condition\/ $C'(\lambda)$}, $\lambda>0$, if
every element of $S$ is a nonempty reduced word, and
for any two distinct words $W_1\kl,W_2\kl\in S$, the length of any
common prefix of $W_1\kl$ and $W_2\kl$ is less than $\lambda\min\{|W_1\kl|,|W_2\kl|\}$.
A symmetrized group presentation $\langle\,\goth A\,\|\,\Cal R\,\rangle$
is said to satisfy the condition $C'(\lambda)$, $\lambda>0$, if
the set of defining words $\Cal R$ satisfies $C'(\lambda)$.
\enddefinition

This and other classical small cancellation conditions and their applications
may be found in \cite{LySc}.

\remark{Remark\/ \rm11.1}
A symmetrized group presentation $\langle\,\goth A\,\|\,\Cal R\,\rangle$
satisfies the condition $C'(\lambda)$, $\lambda>0$,
if and only if for every reduced simple disc diagram $\Delta$
over $\langle\,\goth A\,\|\,\Cal R\,\rangle$,
there exists $\lambda_2\kl<\lambda$
such that $\Delta$ satisfies the condition
$\Cal A(0,0,0;0,\lambda_2\kl,0,0)$
(all the constants here are regarded as constant functions on $\Delta(2)$).
\endremark

\head
12. Proof of Theorem 1
\endhead

\demo{Proof of Theorem\/ {\rm1 (}see Introduction\/{\rm)}\rm}
Let $F_n\kl$ be a free group of rank $n$, $n\in\Bbb N$.
Let $\goth B=\{x_1\kl,\,\ldots,x_n\kl\}$ be a basis of~$F_n\kl$.
Let $\lambda$ be a positive number less than~$1/13$.
Let $S$ be an infinite subset of $F_n\kl$ such that all the
elements of $S$ are cyclically reduced relative to $\goth B$
and are not proper powers in~$F_n\kl$.
Suppose that the symmetrization of $S$ relative to $\goth B$
satisfies the small cancellation condition $C'(\lambda)$ relative to~$\goth B$.
The goal is to prove that there exists an infinite simple $2$-generated group $P$ and
a homomorphism $\phi\:F_n\kl\to P$ such that $\phi$ maps $S$ surjectively onto~$P$.

Let $S'$ be an infinite subset of $S$ which satisfies the following three properties
if the elements of $F_n\kl$ are regarded as reduced group words in the alphabet~$\goth B$:
\medskip

\itemitem{$(1)$}
If $w_1\kl sv_1\kl$ or $(w_1\kl sv_1\kl)^{-1}$ belongs to $S'$, and
$w_2\kl sv_2\kl$ or $(w_2\kl sv_2\kl)^{-1}$ belongs to $S'$, then
either $w_1\kl=w_2\kl$ and $v_1\kl=v_2\kl$,
or $|s|\le\lambda\min\{|w_1\kl sv_1\kl|,|w_2\kl sv_2\kl|\}$.
(Here, the phrase ``$w_1\kl sv_1\kl$ belongs to $S'$'' means
that the concatenation of the group words $w_1\kl$, $s$, $v_1\kl$
is a reduced group word representing an element of $S'$.)
\smallskip

\itemitem{$(2)$}
If $u_1\kl$ and $u_2\kl$ are elements of $S'$, then $u_1\Kl\ne u_2^{-1}$.
\smallskip

\itemitem{$(3)$}
If $u$ is an element of $S'$, then $|u|\ge5$.
\medskip

\noindent
Such a set $S'$ may be obtained from $S$ by first leaving out all the elements
of length less than $5$ and then picking one representative out
of each of the equivalence classes of the following equivalence relation:
call elements $u_1\kl$ and $u_2\kl$ equivalent if
$u_2\Kl$ is a cyclic shift of $u_1\Kl$ or~$u_1^{-1}$.
Note that condition $(1)$ implies that if $u\in S'$,
then any common subword of $u$ and $u^{-1}$ is of length at most $\lambda|u|$
(though, it also follows from the condition $C'(\lambda)$ for $S$).

Let
$$
\lambda_1\kl=\frac{5}{13}-5\lambda,\quad \lambda_2\kl=\lambda,\quad \mu=\frac{5}{13}.
$$
Then $\lambda_1\kl,\lambda_2\kl\in[0,1]$ and
$$
\gathered
2\lambda_1\kl+13\lambda_2\kl<1,\\
\lambda_1\kl+5\lambda_2\kl=\mu.
\endgathered
$$

Let $F_{n+2}\kl$ be a free group of rank $n+2$ with a basis
$\goth A=\{x_1\kl,\,\ldots,x_n\kl,a,b\}$, containing $F_n\kl$ as a subgroup.
Note that $F_{n+2}\kl$ is the inner free product of the subgroup~$F_n\kl$
and the subgroup generated by $\{a,b\}$.
In the rest of this proof, all the elements of $F_{n+2}\kl$
shall be regarded as reduced group words in the alphabet~$\goth A$.
The product of two elements of $F_{n+2}\kl$ shall mean the concatenation
of the group words; it may be not reduced.
The usual group multiplication will not be used.

Split the set $S'$ into the disjoint union of two infinite sets $S_1\kl$ and~$S_2\kl$.
The first step of ``constructing'' the group $P$ consists in
imposing three systems of relations on the group~$F_{n+2}\kl$.
The first system of relations uses elements of $S_1\kl$ and
ensures that every element of the obtained quotient group $G$
is represented by some element of $S_1\kl$
(in particular, the quotient homomorphism maps $F_n\kl$ onto $G$);
the second system uses elements of $S_2\kl$ and ensures
that the quotient group $G$ has no nontrivial finite homomorphic images;
the third system ensures that the elements of $G$
represented by $a$ and $b$ generate the whole of~$G$.
The main theorem shall be used for proving that the quotient group $G$
is not trivial.
The second step of constructing $P$ consists in taking the quotient
of $G$ over its maximal proper normal subgroup.

Let $v_{1,1}\kl$, $v_{1,2}\kl$, $v_{1,3}\kl$, \dots be a list of all reduced group words
in the alphabet~$\goth A$.
Using infiniteness of $S_1\kl$ and finiteness of $\goth B$,
it is easy to show that there exists a system $\{u_{1,i}\kl\}_{i=1}^{+\infty}$
of pairwise distinct elements of $S_1\kl$ such that
$\lambda_1\kl|u_{1,i}\kl|\ge|v_{1,i}\kl|$ for every $i\in\Bbb N$.
For every $i\in\Bbb N$, let
$$
r_{1,i}\Kl=u_{1,i}\Kl v_{1,i}^{-1}.
$$
(Here, $r_{1,i}\kl$ is the concatenation of $u_{1,i}\kl$ and $v_{1,i}^{-1}$;
it does not need to be reduced.)
Let $\Cal R_1\kl=\{\,r_{1,i}\kl\mid i=1,2,\ldots\,\}$.

Let a function from a set $X$ to a group $G$ be called {\it trivial\/} if
it maps all the elements of $X$ to the neutral element of~$G$.
Let $M$ be a set of nontrivial finite groups such that
every nontrivial finite group is isomorphic to exactly one group in~$M$.
Such a set $M$ exists and is countable.
Informally speaking, $M$ is a set of up to isomorphism all nontrivial finite groups.
Let $T$ be the set of all ordered pairs ${(}G,\psi{)}$
where $G\in M$, and $\psi$ is a nontrivial function from $\goth A$ to~$G$.
Clearly, $T$ is countable.
Let ${(} G_1\kl,\psi_1\kl{)}$, ${(} G_2\kl,\psi_2\kl{)}$, \dots
be a list of all elements of $T$ (without recurrences).
Let $u_{2,1}\kl$, $u_{2,2}\kl$, $u_{2,3}\kl$, \dots be a list (without recurrences)
of elements of $S_2\kl$ of length at least $1/\lambda_1\kl$.
For every $i\in\Bbb N$, let $v_{2,i}\kl$ be a group word over $\goth A$ of minimal length
such that the values of $v_{2,i}\kl$ and $u_{2,i}\kl$ in $G_i\kl$ with respect to $\psi_i\kl$
are not equal.
Clearly, the length of every such $v_{2,i}\kl$ is either $0$ or~$1$.
For every $i\in\Bbb N$, let
$$
r_{2,i}\Kl=u_{2,i}\Kl v_{2,i}^{-1}.
$$
Let $\Cal R_2\kl=\{\,r_{2,i}\kl\mid i=1,2,\ldots\,\}$.

Let $S_3\kl=\{u_{3,1}\kl,\,\ldots,u_{3,n}\kl\}$ be a set
consisting of $n$ reduced group words
of length at least $1/\lambda_1\kl$ in the alphabet $\{a,b\}$
which satisfies the same three conditions as those required of $S'$ above.
(Such a set $S_3\kl$ exists.)
For every $i\in\{1,\,\ldots,n\}$, let $v_{3,i}\kl=x_i\kl$ and
$$
r_{3,i}\Kl=u_{3,i}\Kl v_{3,i}^{-1}.
$$
Let $\Cal R_3\kl=\{\,r_{3,i}\kl\mid i=1,\,\ldots,n\,\}$.

Let
$$
\Cal R=\Cal R_1\kl\cup\Cal R_2\kl\cup\Cal R_3\kl.
$$
Let $G$ be the group defined by the presentation $\langle\,\goth A\,\|\,\Cal R\,\rangle$.
This group may be naturally identified with a quotient group of~$F_{n+2}\kl$.
Let $\phi_1\kl\:F_{n+2}\kl\to G$ be the corresponding epimorphism.

Because of the relations $r=1$, $r\in\Cal R_1\kl$,
the homomorphism $\phi_1\kl$ maps $S_1\kl$ onto~$G$.

Because of the relations $r=1$, $r\in\Cal R_2\kl$,
the group $G$ has no nontrivial finite quotients.
Indeed, suppose $G$ has a nontrivial finite quotient.
Then for some $i$ there exists an epimorphism $\phi_2\kl$
from $G$ onto $G_i\kl$ such that $\phi_2\kl\circ\phi_1\kl$ extends~$\psi_i\kl$.
On one hand, the value of the word $r_{2,i}\kl$ in $G_i\kl$ with respect to $\psi_i\kl$
is nontrivial by the construction of~$r_{2,i}\kl$.
On the other hand, the value of $r_{2,i}\kl$ in $G$ with respect to $\phi_1\kl$ is $1$
since $r_{2,i}\kl\in\Cal R$;
therefore, the value of $r_{2,i}\kl$ in $G_i\kl$ with respect to $\phi_2\kl\circ\phi_1\kl$ is~$1$.
This gives a contradiction.

Because of the relations $r=1$, $r\in\Cal R_3\kl$,
the group $G$ is generated by $\phi_1\kl(a)$ and $\phi_1\kl(b)$.

The only property of $G$ that still needs to be established,
is that $G$ is nontrivial.

Suppose $G$ is trivial.
Then, by van Kampen's Lemma, for any group word $v$ in the alphabet $\goth A$,
there exists a disc diagram $\Delta$ over the presentation
$\langle\,\goth A\,\|\,\Cal R\,\rangle$ such that the label of some
representative of the contour of $\Delta$ is~$v$.
In particular, there exists a reduced disc diagram over
$\langle\,\goth A\,\|\,\Cal R\,\rangle$ whose contour has length~$1$.
Let $\Delta$ by such a diagram.
Clearly, it is simple.

Define a selection on $\Delta$ as follows.
Consider all the faces of $\Delta$ one-by-one.
Take an arbitrary face $f$ of $\Delta$ on which the selection
has not been defined yet.
Let $i$ and $j$ be such indices that the label of
some representative of the contour of $f$ is
either $r_{i,j}\Kl$ or~$r_{i,j}^{-1}$.
Then there exist paths $p$ and $q$ such that $pq$ is a representative
of the contour of $f$, and
either the label of $p$ is $u_{i,j}\kl$, the label of $q$ is $v_{i,j}^{-1}$,
or the label of $p$ is $u_{i,j}^{-1}$, the label of $q$ is~$v_{i,j}\kl$.
In both cases, such a path $p$ is reduced since the word $u_{i,j}\kl$ is reduced.
Let the selection on $f$ consist of all the nontrivial subpaths of this~$p$.
The diagram $\Delta$ has become a {\it diagram with selection\/}
or an {\it S-diagram}.

Now, check that the S-diagram $\Delta$ satisfies the condition
$\Cal A(1,1,1;\lambda_1\kl,\lambda_2\kl,0,0)$
(here, all the constants shall be interpreted as constant functions on $\Delta(2)$).

The condition $\Cal A_1\kl(1)$ follows directly from the construction of the selection.
The conditions $\Cal A_2\kl(1)$ and $\Cal A_3\kl(1)$ follow from $\Cal A_1\kl(1)$
(see Remark 4.1).

In order to check the condition $\Cal A_4\kl(\lambda_1\kl)$,
consider an arbitrary face $f$ and the corresponding subpaths $p$ and $q$ of its contour
used above in defining the selection on~$f$.
Note that $p$ is the only maximal selected subpath of the contour of $f$,
and $q$ is the ``complement'' of $p$ in the contour of~$f$.
Take $i$ and $j$ such that $p$ is labelled by $u_{i,j}\Kl$ or $u_{i,j}^{-1}$,
and $q$ is labelled by $v_{i,j}^{-1}$ or~$v_{i,j}\Kl$.
By the construction of the words $u_{i,j}\kl$ and $v_{i,j}\kl$, the inequality
$\lambda_1\kl|u_{i,j}\kl|\ge|v_{i,j}\kl|$ holds.
Therefore,
$$
|q|\le\lambda_1\kl|p|\le\lambda_1\kl|pq|=\lambda_1\kl d(f).
$$
The condition $\Cal A_4\kl(\lambda_1\kl)$ follows.

To check the conditions
$\Cal A_5\kl(\lambda_2\kl)$, $\Cal A_6\kl(\lambda_2\kl,0)$, $\Cal A_7\kl(\lambda_2\kl,0)$,
it is enough to prove that
for all faces $f$ of $\Delta$, the length of every double-selected
oriented arc $t$ which is a subpath of the contour of $f$
is at most~$\lambda_2\kl d(f)$.

Consider first an arbitrary intro-facial double-selected oriented arc~$t$.
Let $f$ be the face whose contour has $t$ as a subpath.
Let $p$ be the maximal selected subpath of the contour of $f$,
and $q$ be such that $pq$ is a representative of the contour of $f$
(the same notation as above).
Let $u$ be the label of~$p$.
By the choice of the selection, either $u$ or $u^{-1}$ belongs to~$S'\cup S_3\kl$.
The contour of $f$ has the form $ts_1\kl t^{-1}s_2\kl$.
Since $t$ and $t^{-1}$ are selected, and there is exactly one
maximal selected subpath of the contour of $f$,
either $ts_1\kl t^{-1}$ or $t^{-1}s_2\kl t$ is selected.
So, either $ts_1\kl t^{-1}$ or $t^{-1}s_2\kl t$ is a subpath of~$p$.
In either case, the label of $t$ is a common subword of $u$ and~$u^{-1}$.
Therefore,
$$
|t|\le\lambda|p|\le\lambda d(f)=\lambda_2\kl d(f)
$$
(see property (1) of $S'$ and $S_3\kl$).

Consider next an arbitrary inter-facial double-selected oriented arc~$t$.
Let $f_1\kl$ be the face whose contour has $t$ as a subpath,
and $f_2\kl$ be the face whose contour has $t^{-1}$ as a subpath
($t$ is between $f_1\kl$ and~$f_2\kl$).
Let $p_1\kl$ be the maximal selected subpath of the contour of $f_1\kl$,
and $p_2\kl$ be the maximal selected subpath of the contour of~$f_2\kl$.
By the choice of the selection, $|p_1\kl|\le d(f_1\kl)$ and $|p_2\kl|\le d(f_2\kl)$.
Let $u_1\kl$ be the label of $p_1\kl$, and $u_2\kl$ be the label of~$p_2\kl$.
One of the words $u_1\Kl$, $u_1^{-1}$ and one of the words $u_2\Kl$, $u_2^{-1}$
belong to~$S'\cup S_3\kl$.
The label of $t$ is a common subword of $u_1\Kl$ and~$u_2^{-1}$.
If $|t|>\lambda|p_1\kl|$ or $|t|>\lambda|p_2\kl|$, then
the pair of faces $\{f_1\kl,f_2\kl\}$ is cancellable
(follows from the properties of $S'$ and $S_3\kl$, the fact that
$S'$ and $S_3\kl$ use disjoint alphabets, and the construction of $\Cal R$).
Since $\Delta$ is reduced, this cannot happen.
Therefore,
$$
|t|\le\lambda|p_1\kl|\le\lambda d(f_1\kl)=\lambda_2\kl d(f_1\kl).
$$
Thus, $\Delta$ satisfies the conditions
$\Cal A_5\kl(\lambda_2\kl)$, $\Cal A_6\kl(\lambda_2\kl,0)$, $\Cal A_7\kl(\lambda_2\kl,0)$.

Since $\Delta$ satisfies the condition $\Cal A(1,1,1;\lambda_1\kl,\lambda_2\kl,0,0)$,
by the main theorem and subsequent Remark 6.3,
the length of its contour is at least
$(1-2\mu)\sum_{f\in\Delta(2)}d(f)$.
Since
$$
\aligned
(1-2\mu)\sum_{f\in\Delta(2)}d(f)
&\ge(1-2\mu)\min\bigl\{\,|r|\bigm|r\in\Cal R\,\bigr\}\\
&\ge(1-2\mu)\min\bigl\{\,|u|\bigm|u\in S_1\kl\cup S_2\kl\cup S_3\kl\,\bigr\}
\ge\frac{3}{13}\cdot5>1,
\endaligned
$$
this contradicts to the assumption that the length of the contour of $\Delta$ is~$1$.
Hence, the group $G$ could not be trivial.

By Zorn's Lemma, there exists a maximal proper normal subgroup $N$ of the group $G$
since $G$ is finitely generated.
Let $P=G/N$.
Let $\phi_2\kl\:G\to P$ be the quotient homomorphism.
Then $P$ is a desired infinite simple $2$-generated group,
and $\phi=\phi_2\kl\circ\phi_1\kl|_{F_n\kl}\Kl$ is a desired epimorphism $F_n\kl\to P$.
\qed
\enddemo

\demo{Proof of the Corollary of Theorem\/ {\rm1 (}see Introduction\/{\rm)}\rm}
Take a free group $F$ of rank $27$ with a basis $\{a_1\kl,\,\ldots,a_{27}\kl\}$.
Let
$$
S=\{\,a_1^m\cdots a_{27}^m\mid m\in\Bbb N\,\}.
$$
The symmetrization of the set $S$ obviously satisfies the condition
$C'(\frac{2}{27}+\varepsilon)$
relative to the basis $\{a_1\kl,\,\ldots,a_{27}\kl\}$ for any $\varepsilon>0$.
Therefore, by Theorem~1, there exists an infinite simple $2$-generated group $G$ and
a homomorphism $\phi\:F\to G$ such that $\phi$ maps $S$ onto~$G$.
Let $x_i\kl=\phi(a_i\kl)$, $i=1$, \dots,~$27$.
\qed
\enddemo

\head
13. Proof of Theorem 2
\endhead

\demo{Proof of Theorem\/ {\rm2 (}see Introduction\/{\rm)}\rm}
In this proof, let
$$
k_1\kl=28,\quad
\lambda_1\kl=10^{-2},\quad \lambda_2\kl=10^{-4},\quad \lambda_4\kl=\frac{1}{28},\quad
\lambda_3\kl=4\lambda_4\kl=\frac{1}{7}.
$$
Let $\mu=0.45$.
These numbers satisfy the inequalities
$$
\gathered
2\lambda_1\kl+(8+5k_1\kl)\lambda_2\kl+6\lambda_3\kl+2\lambda_4\kl<1,\\
\lambda_1\kl+(3+2k_1\kl)\lambda_2\kl+3\lambda_3\kl\le\mu.
\endgathered
$$
Let $\bar a={(} k_1\kl,k_1\kl,k_1\kl;\lambda_1\kl,\lambda_2\kl,\lambda_3\kl,\lambda_4\kl{)}$.

Let $\goth A=\{x,y\}$ be a $2$-letter alphabet.
Let ${(} a_1\kl,b_1\kl{)}$, ${(} a_2\kl,b_2\kl{)}$, \dots be all the ordered pairs
of reduced group words in the alphabet $\goth A$
listed without recurrences.
Let, moreover, the function $i\mapsto{(} a_i\kl,b_i\kl{)}$ be computable.

Pick two nonempty reduced words $v_1\kl$, $v_2\kl$ and a system of nonempty
reduced words $\{u_{i,j}\kl\}_{i=1,2,\ldots\hfill\atop j=1,\,\ldots,14}\kl$
so that the following conditions are satisfied:
\medskip

\itemitem{$(1)$}
If $i\in\Bbb N$, then
$|u_{i,1}\kl|=|u_{i,2}\kl|=\cdots=|u_{i,14}\kl|$.
\smallskip

\itemitem{$(2)$}
If $i\in\Bbb N$, then
$28|u_{i+1,1}\kl|+14|a_{i+1}\kl|+|b_{i+1}\kl|>28|u_{i,1}\kl|+14|a_i\kl|+|b_i\kl|$.
\smallskip

\itemitem{$(3)$}
If $i\in\Bbb N$, then
$28|u_{i,1}\kl|\ge99(14|a_i\kl|+|b_i\kl|)$.
\smallskip

\itemitem{$(4)$}
If $i_1\kl,i_2\kl\in{\Bbb N}$, $j_1\kl,j_2\kl\in\{1,2,3,\,\ldots,14\}$,
$\sigma_1\kl,\sigma_2\kl\in\{\pm1\}$,
$u_{i_1\kl,j_1\kl}^{\sigma_1\kl}=w_1\Kl sz_1\Kl$ and
$u_{i_2\kl,j_2\kl}^{\sigma_2\kl}=w_2\Kl sz_2\Kl$, then
either $w_1\kl=w_2\kl$, $z_1\kl=z_2\kl$, $i_1\kl=i_2\kl$, $j_1\kl=j_2\kl$, $\sigma_1\kl=\sigma_2\kl$,
or $|s|\le28\cdot10^{-4}\min\{|u_{i_1\kl,j_1\kl}\kl|,|u_{i_2\kl,j_2\kl}\kl|\}$.
\smallskip

\itemitem{$(5)$}
The products (\ie, concatenations)
$v_1\kl v_1\kl$, $v_2\kl v_2\kl$, $v_1\kl v_2\kl$,
$v_1^{-1}v_2\Kl$, $v_1\Kl v_2^{-1}$, $v_1^{-1}v_2^{-1}$
are reduced.
For example, $v_1\kl$ starts and ends with the letter $x$, and
$v_2\kl$ starts and ends with the letter~$y$.
\smallskip

\itemitem{$(6)$}
If $v\in\{v_1^{\pm1},v_2^{\pm1}\}$ and
$u\in\{\,u_{i,j}^{\pm1}\mid i\in{\Bbb N},j=1,\,\ldots,14\,\}$, then
$|u|\ge40|v|$ and $v$ is not a subword of~$u$.
\medskip

\noindent
Besides that, assume that the function ${(} i,j{)}\mapsto u_{i,j}\kl$
is computable.
Such words $v_1\kl$, $v_2\kl$ and a system
$\{u_{i,j}\kl\}_{i=1,2,\ldots\hfill\atop j=1,\,\ldots,14}\kl$ exist.

For every $i=1$, $2$, $3$, \dots, let
$$
r_i\Kl=a_i^{u_{i,1}\kl}a_i^{u_{i,2}\kl}\cdots a_i^{u_{i,14}\kl}b_i^{-1}
$$
(recall that $a_i^{u_{i,j}\kl}=u_{i,j}\Kl a_i\Kl u_{i,j}^{-1}$).
Note that by condition $(2)$,
$$
|r_1\kl|<|r_2\kl|<|r_3\kl|<\cdots,
$$
by condition $(3)$, for every $i$,
$$
14|a_i\kl|+|b_i\kl|\le\lambda_1\kl|r_i\kl|,
$$
by condition $(1)$, for every $i$,
$$
|u_i\kl|\le\lambda_4\kl|r_i\kl|.
$$

Define a sequence of sets $\{\Cal R_i\kl\}_{i=0}^{+\infty}$ inductively.
Let $\Cal R_0\kl=\emptyset$.
Take $i>0$.
If the relation $a_i\kl=1$ is a consequence of the relations $r=1$, $r\in\Cal R_{i-1}\kl$
(for example, if $a_i\kl$ is the empty word), then let $\Cal R_i\kl=\Cal R_{i-1}\kl$.
Otherwise, let $\Cal R_i\kl=\Cal R_{i-1}\kl\cup\{r_i\kl\}$.
Let
$$
\Cal R=\bigcup_{i=1}^{+\infty}\Cal R_i\kl.
$$

Let $G$ be the group defined by the presentation $\langle\,\goth A\,\|\,\Cal R\,\rangle$.
Let $\phi$ be the natural function from the set of all group words
in the alphabet $\goth A$ onto~$G$.
Let $H$ be the subgroup of $G$ generated by $\phi(v_1\kl)$ and $\phi(v_2\kl)$.
It easily follows from the imposed relations that $G$ is $14$-boundedly simple.
Moreover, for any nontrivial element $g\in G$ and for any element $h\in G$,
the element $h$ is the product of $14$ conjugates of~$g$.
The next goal is to prove that the word problem in $G$ is solvable and that
$H$ is freely generated by $\phi(v_1\kl)$ and $\phi(v_2\kl)$.
It will be done using the fact that every reduced simple disc diagram $\Delta$
over $\langle\,\goth A\,\|\,\Cal R\,\rangle$ satisfies the condition
$\Cal A(\bar a)$
(where the components of $\bar a$ are viewed as constant functions on $\Delta(2)$).

In the rest of this proof, a selection on a diagram $\Delta$
over $\langle\,\goth A\,\|\,\Cal R\,\rangle$ will be called {\it special\/} if
for every face $f$ of $\Delta$,
there is a natural number $i$, a representative $c$ of the contour of $f$,
and paths $q_0\kl$, \dots, $q_{14}\kl$,
$p_1\kl$, \dots, $p_{14}\kl$, $p_1'$, \dots, $p_{14}'$
such that the following conditions hold:
\medskip

\item{$\bullet$}
the label of $c$ is either $r_i\Kl$ or~$r_i^{-1}$;

\item{$\bullet$}
if the label of $c$ is $r_i\kl$, then
$$
c=p_1'q_1\Kl p_1^{-1}p_2'q_2\Kl p_2^{-1}\cdots
p_{14}'q_{14}\Kl p_{14}^{-1}q_0^{-1},
$$
and the selection on $f$ consists of all nontrivial subpaths
of $p_1^{-1}$, \dots, $p_{14}^{-1}$ and $p_1'$, \dots,~$p_{14}'$;

\item{$\bullet$}
if the label of $c$ is $r_i^{-1}$, then
$$
c=q_0\Kl p_{14}\Kl q_{14}^{-1}p_{14}^{\prime-1}p_{13}\Kl q_{13}^{-1}p_{13}^{\prime-1}\cdots
p_1\Kl q_1^{-1}p_1^{\prime-1},
$$
and the selection on $f$ consists of all nontrivial subpaths
of $p_1^{\prime-1}$, \dots, $p_{14}^{\prime-1}$ and $p_1\kl$, \dots,~$p_{14}\kl$.

\item{$\bullet$}
the label of $q_0\kl$ is $b_i\kl$, and for every $j=1$, \dots, $14$,
the labels of $p_j\kl$ and $p_j'$ are $u_{i,j}\kl$, the label of $q_j\kl$ is~$a_i\kl$.
\medskip

\noindent
Clearly, on every diagram over $\langle\,\goth A\,\|\,\Cal R\,\rangle$
there exists a special selection.
A diagram over $\langle\,\goth A\,\|\,\Cal R\,\rangle$
together with a special selection will be called a {\it special S-diagram}.
Note that if a simple disc S-map $\Delta'$ is cut out from a special S-diagram,
and the S-map $\Delta'$ is naturally labelled via some pasting morphism,
then the S-diagram obtained from $\Delta'$ is special.

Consider an arbitrary special S-diagram~$\Delta$.
Clearly, $\Delta$ satisfies the conditions
$\Cal A_1\kl(k_1\kl)$, $\Cal A_2\kl(k_1\kl)$, and $\Cal A_3\kl(k_1\kl)$.
It easily follows from condition $(3)$ that $\Delta$ also satisfies
$\Cal A_4\kl(\lambda_1\kl)$.
It is to be shown that in fact every reduced special disc S-diagram
satisfies $\Cal A(\bar a)$.

Use induction on the number of internal edges of a reduced special disc S-diagram.
If a reduced special S-diagram has no internal edges, then it clearly satisfies
the condition $\Cal A(\bar a)$.
Suppose that $\Delta$ is a reduced special disc S-diagram,
and every reduced special disc S-diagram having fewer internal edges than $\Delta$
satisfies the condition $\Cal A(\bar a)$.
To complete the inductive proof, it is enough to prove now that $\Delta$
satisfies $\Cal A(\bar a)$.
Since it is enough to prove that every maximal simple disc S-subdiagram
of $\Delta$ satisfies $\Cal A(\bar a)$, assume that $\Delta$ is simple.

The length of every selected subpath of the contour of every face $f$
of $\Delta$ is at most $\lambda_4\kl d(f)$.
Consider arbitrary (not necessarily distinct) faces $f_1\kl$ and $f_2\kl$ of~$\Delta$.
Suppose there exists a double-selected arc $t$ between $f_1\kl$ and $f_2\kl$
such that $|t|>\lambda_2\kl d(f_1\kl)$
(call such arcs ``long'').
Then, by condition $(4)$, the labels of $f_1\kl$ and $f_2\kl$ are either equal
or mutually inverse.
Since the diagram $\Delta$ is reduced, the labels of $f_1\kl$ and $f_2\kl$ in fact are equal.
Hence, the contours of $f_1\kl$ and $f_2\kl$ have the same length.
Therefore, if $p$ is a selected subpath of the contour of $f_2\kl$
on which a given ``long'' arc lies,
then the sum of the lengths of all maximal double-selected arcs
between $f_1\kl$ and $f_2\kl$ that lie on $p$ is at most $\lambda_4\kl d(f_1\kl)$.
Suppose, moreover, that the faces $f_1\kl$ and $f_2\kl$ are distinct.
Then it follows from the choice of the selection and from the properties of the system
$\{u_{i,j}\kl\}_{i=1,2,\ldots\hfill\atop j=1,\,\ldots,14}\kl$
that all such ``long'' double-selected arcs between $f_1\kl$ and $f_2\kl$
cannot be situated on more
than $4$ maximal selected subpaths of the contour of $f_1\kl$ (as well as of~$f_2\kl$).
Therefore, the sum of the lengths of all such ``long''  maximal double-selected arcs
between $f_1\kl$ and $f_2\kl$ is not greater than $4\lambda_4\kl d(f_1\kl)=\lambda_3\kl d(f_1\kl)$.
Thus, $\Delta$ satisfies the conditions $\Cal A_7\kl(\lambda_2\kl,\lambda_4\kl)$
and $\Cal A_6\kl(\lambda_2\kl,\lambda_3\kl)$.

Verification of the condition $\Cal A_5\kl(\lambda_2\kl)$ is more complicated.
Suppose that $\Cal A_5\kl(\lambda_2\kl)$ does not hold in~$\Delta$.
Let $f_0\kl$ be a face of $\Delta$ and $t$ be a double-selected inward
intro-facial oriented arc such that
$t$ is a subpath of the contour of $f_0\kl$ and $|t|>\lambda_2\kl d(f_0\kl)$
(such $f_0\kl$ and $t$ exist since $\Cal A_5\kl(\lambda_2\kl)$ does not hold).
Let $s_1\kl$ and $s_2\kl$ be such subpaths of the contour of $f_0\kl$ that
$ts_1\kl t^{-1}s_2\kl$ is a representative of the contour of~$f_0\kl$.
Since every reduced special disc S-diagram having fewer internal edges than $\Delta$
satisfies the condition $\Cal A_5\kl(\lambda_2\kl)$,
the choice of $f_0\kl$ and $t$ is in fact unique
and the path $s_2^{-1}$ is a representative of the contour of $\Delta$
(otherwise the cycle represented by $s_2^{-1}$ would cut out
an S-diagram having fewer internal edges and still
not satisfying $\Cal A_5\kl(\lambda_2\kl)$).

Let $i$ be the number such that $d(f_0\kl)=|r_i\kl|$.
Then the label of some representative of the contour of $f_0\kl$
is either $r_i\kl$ or~$r_i^{-1}$.
By the construction of $r_i\kl$ (see condition $(4)$)
and by the definition of a special selection,
there is exactly one $m\in\{1,2,\,\ldots,14\}$ and one $\sigma\in\{\pm1\}$
such that the label of $t$ is a subword of~$u_{i,m}^\sigma$.
Then either the label of $s_1\kl$ or the label of $t^{-1}s_2\kl t$ has the form
$$
\bigl(wa_i\kl w^{-1}\bigr)^{\pm1},
$$
where $w$ is a suffix of~$u_{i,m}\kl$.

Case 1: the label of $s_1\kl$ has the form
$$
\bigl(wa_i\kl w^{-1}\bigr)^{\pm1},
$$
$w$ is a suffix of~$u_{i,m}\kl$.
Let $\Delta_0\kl$ be the disc subdiagram of $\Delta$
with the contour represented by~$s_1\kl$.
Note that $\Delta(2)=\Delta_0\kl(2)\cup\{f_0\kl\}$.
By the construction of $\{\Cal R_i\kl\}_{i=0}^{+\infty}$,
the label of $s_1\kl$ is not trivial modulo the set of relation
$$
\{\,r=1\mid r\in\Cal R_{i-1}\kl\,\}=\{\,r=1\mid r\in\Cal R,\ |r|<d(f_0\kl)\,\}.
$$
Pick a face $f_1\kl$ in $\Delta_0\kl$ such that $d(f_1\kl)\ge d(f_0\kl)$.
Let $\Delta_1\kl$ be the maximal simple disc S-subdiagram of $\Delta$
containing the face $f_1\kl$ and not containing the face~$f_0\kl$.
All the external edges of $\Delta_1\kl$ lie on the path~$s_1\kl$.
In particular, the length of the contour of $\Delta_1\kl$ is at most~$|s_1\kl|$.
The S-map $\Delta_1\kl$ satisfies $\Cal A_5\kl(\lambda_2\kl)$ since it has
fewer internal edges than~$\Delta$.
Hence, by the main theorem and subsequent Remark 6.3, the length of the contour of $\Delta_1\kl$
is at least $(1-2\mu)d(f_1\kl)$.
Therefore, on one hand,
$$
|s_1\kl|\ge(1-2\mu)d(f_1\kl)\ge(1-2\mu)d(f_0\kl)=\frac{1}{10}d(f_0\kl).
$$
On the other hand,
$$
|s_1\kl|<|a_i\kl|+2|u_{i,m}\kl|\le\frac{1}{14}d(f_0\kl).
$$
Contradiction.

Case 2: the label of $t^{-1}s_2\kl t$ has the form
$$
\bigl(wa_i\kl w^{-1}\bigr)^{\pm1},
$$
$w$ is a suffix of~$u_{i,m}\kl$.
Let $\Delta'$ be an S-diagram cut out of $\Delta$
by the cycle represented by $t^{-1}s_2^{-1}t$.
(Informally speaking, the S-diagram $\Delta'$ may be thought of as
obtained from $\Delta$ by cutting $\Delta$ along $t$
from the boundary inside.)
The S-diagram $\Delta'$ is a simple reduced special disc S-diagram,
it has fewer internal edges than~$\Delta$.
Therefore, it satisfies $\Cal A_5\kl(\lambda_2\kl)$.
Hence, by the main theorem and subsequent Remark 6.3,
the length of the contour of $\Delta'$
is at least $(1-2\mu)d(f_0')$ where $f_0'$ is
the face of $\Delta'$ corresponding to $f_0\kl$ (note that $d(f_0')=d(f_0\Kl)$).
Therefore, on one hand,
$$
|t^{-1}s_2\kl t|\ge(1-2\mu)d(f_0\kl)=\frac{1}{10}d(f_0\kl).
$$
On the other hand,
$$
|t^{-1}s_2\kl t|\le|a_i\kl|+2|u_{i,m}\kl|\le\frac{1}{14}d(f_0\kl).
$$
Contradiction.

Thus, every reduced special disc S-diagram $\Delta$
satisfies $\Cal A(\bar a)$.
In particular, $\Delta$, $\bar a$, and $\mu$ satisfy the
assumptions of the main theorem.
Hence, if $n$ is the length of the contour of $\Delta$, then
the sum of the lengths of the contours of all the faces of $\Delta$
is at most $\frac{1}{1-2\mu}n$, and the number of edges of $\Delta$
is at most $\frac{1-\mu}{1-2\mu}n$.

Now, if $\Delta$ is an arbitrary reduced disc diagram over
$\langle\,\goth A\,\|\,\Cal R\,\rangle$,
the length of the contour of $\Delta$ is greater than or equal to
the sum of the lengths of the contours of all its maximal simple disc subdiagrams.
Hence, if $n$ is the length of the contour of $\Delta$, then again
the sum of the lengths of the contours of all the faces of $\Delta$
is at most $\frac{1}{1-2\mu}n$, and the number of edges of $\Delta$
is at most $\frac{1-\mu}{1-2\mu}n$.

For every $i=0$, $1$, $2$, \dots, let $G_i\kl$ be the group defined by the presentation
$\langle\,\goth A\,\|\,\Cal R_i\kl\,\rangle$,
and let $\phi_i\kl$ be the natural function from the set of all group words
in the alphabet $\goth A$ onto~$G_i\kl$.

Consider an arbitrary group word $w$ whose value in $G$ is~$1$.
Let $\Delta$ be a reduced deduction diagram for $w$
over $\langle\,\goth A\,\|\,\Cal R\,\rangle$.
Let $n=|w|$.
Consider an arbitrary $i>10n$.
Then $|r_i\kl|>i>10n=(1-2\mu)^{-1}n$.
So, no face in $\Delta$ is labelled with a cyclic word associated with
$r_i\kl$ or~$r_i^{-1}$.
Hence, the value of $w$ in $G_{10n}\kl$ is $1$ as well.

Recall that an algorithm solves the word problem in
$\langle\,\goth A\,\|\,\Cal R\,\rangle$ if
for every group word $w$ in the alphabet $\goth A$,
it determines whether or not $w$ is trivial in $G$,
\ie, whether the relation $w=1$ is a consequence
of the relations $r=1$, $r\in\Cal R$.
Here is an informal description of an algorithm that solves the word
problem in $\langle\,\goth A\,\|\,\Cal R\,\rangle$.

Consider group words only in the alphabet~$\goth A$.
Let {\tt Alg\_1} be an algorithm that for every finite set $\Cal S$
of group words and for every group word $w$ decides whether
there exists a deduction diagram $\Delta$ for $w$
over $\langle\,\goth A\,\|\,\Cal S\,\rangle$ with at most
$\frac{1-\mu}{1-2\mu}|w|$ edges.
Clearly, such an algorithm exists.
Let {\tt Alg\_2} be an algorithm that for every natural $i$
constructs the set~$\Cal R_i\kl$.
It starts with $\Cal R_0\kl=\emptyset$, and then for every $l=1$, \dots, $i$
constructs $\Cal R_l\kl$ by using {\tt Alg\_1} to check
whether $\Cal R_l\kl=\Cal R_{l-1}\kl$ or $\Cal R_l\kl=\Cal R_{l-1}\kl\cup\{r_l\kl\}$.
Let {\tt Alg\_3} be an algorithm that for every group word $w$
first uses {\tt Alg\_2} to construct $\Cal R_{10n}\kl$ where $n=|w|$,
and then uses {\tt Alg\_1} to determine whether
there exists a deduction diagram $\Delta$ for $w$
over $\langle\,\goth A\,\|\,\Cal R_{10n}\kl\,\rangle$ with at most
$\frac{1-\mu}{1-2\mu}n$ edges.
Such a diagram exists if and only if $w$ is trivial in~$G_{10n}\kl$.
Hence, the algorithm {\tt Alg\_3} solves the word problem
in $\langle\,\goth A\,\|\,\Cal R\,\rangle$.

It is left to prove that $H$ is freely generated by $\phi(v_1\kl)$ and $\phi(v_2\kl)$.
Suppose it is not.
Then there exist
$m\ge1$, $i_1\kl,\,\ldots,i_m\kl\in\{1,2\}$, $\sigma_1\kl,\,\ldots,\sigma_m\kl\in\{\pm1\}$ such that
$$
\bigl(\forall j\in\{1,\,\ldots,m-1\}\bigr)
\bigl(i_j\kl\ne i_{j+1}\kl\text{ or }\sigma_j\kl=\sigma_{j+1}\kl\bigr)
\text{ and }\bigl(i_m\kl\ne i_1\kl\text{ or }\sigma_m\kl=\sigma_1\kl\bigr),
$$
and the group word $w=v_{i_1\kl}^{\sigma_1\kl}\cdots v_{i_m\kl}^{\sigma_m\kl}$
is trivial in~$G$.
By the choice of $v_1\kl$ and $v_2\kl$, the group word $w$ is cyclically reduced.
Consider an arbitrary reduced deduction diagram $\Delta$ for $w$
over $\langle\,\goth A\,\|\,\Cal R\,\rangle$.
The contour of $\Delta$ is reduced since $w$ is cyclically reduced.
Let $\Delta_0\kl$ be a maximal simple disc submap of $\Delta$
such that some representative of the contour $\Delta_0\kl$ is a subpath
of the contour of~$\Delta$.
Let $q_0\kl$ be a representative of the contour of $\Delta_0\kl$
which is a subpath of the contour of~$\Delta$.
Define a selection on $\Delta_0\kl$ the same way as above (special selection).
As it has been shown, the obtained S-diagram satisfies $\Cal A(\bar a)$.
Apply the main theorem with subsequent Remark 6.3 and
the lemma about exposed face to~$\Delta_0\kl$.
Let $f$ be a face of $\Delta_0\kl$ and $P$ be a non-overlapping set
of selected external arcs of $\Delta_0\kl$ incident to $f$ such that
$\|P\|\le k_1\kl+1$ and
$$
\sum_{p\in P}|p|
\ge\bigl(1-2\mu-(2+k_1\kl)\lambda_2\kl-\lambda_4\kl+\|P\|\lambda_2\kl\bigr)d(f).
$$
Note that elements of $P$ are arcs in $\Delta_0\kl$
but do not need to be arcs in $\Delta$
since some intermediate vertex of some element of $P$
may have degree greater than $2$ in~$\Delta$.
Let $P_0\kl$ be the set of all maximal elements of the set
of all arcs that lie on $q_0\kl$ and are subarcs of elements of~$P$.
Then $P_0\kl$ is a non-overlapping set of arcs,
$\|P_0\kl\|\le\|P\|+1$, and $\sum_{p\in P_0\kl}|p|=\sum_{p\in P}|p|$.
Let $p$ be a longest arc in~$P_0\kl$.
Then
$$
\aligned
|p|
&\ge\frac{1-2\mu-(2+k_1\kl)\lambda_2\kl-\lambda_4\kl+\|P\|\lambda_2\kl}{\|P_0\kl\|}d(f)\\
&\ge\frac{1-2\mu-(2+k_1\kl)\lambda_2\kl-\lambda_4\kl+(\|P_0\kl\|-1)\lambda_2\kl}{\|P_0\kl\|}d(f)\\
&=\frac{1-2\mu-(3+k_1\kl)\lambda_2\kl-\lambda_4\kl}{\|P_0\kl\|}d(f)+\lambda_2\kl d(f).
\endaligned
$$
The picked values of $k_1\kl$, $\lambda_2\kl$, $\lambda_4\kl$ and $\mu$
satisfy the inequality
$$
1-2\mu-(3+k_1\kl)\lambda_2\kl-\lambda_4\kl\ge0.
$$
Therefore,
$$
\aligned
|p|
&\ge\frac{1-2\mu-(3+k_1\kl)\lambda_2\kl-\lambda_4\kl}{k_1\kl+2}d(f)+\lambda_2\kl d(f)\\
&=\frac{1-2\mu-\lambda_2\kl-\lambda_4\kl}{k_1\kl+2}d(f).
\endaligned
$$
Therefore, $|p|>\frac{1}{500}d(f)$.
Let $p'$ be the subpath of the contour of $f$ corresponding
to the arc $p$ (so, $p'$ is a selected oriented arc of~$\Delta$).
Let $i$ be the number such that $d(f)=|r_i\kl|$.
Take $j\in\{1,2,\,\ldots,14\}$ and $\sigma\in\{\pm1\}$ such that
the label of $p'$ is a subword of~$u_{i,j}^\sigma$.
Since $28|u_{i,j}\kl|<d(f)$, it follows that
$$
|p'|>\frac{1}{20}|u_{i,j}\kl|\ge2\max\{|v_1\kl|,|v_2\kl|\}.
$$
Since $p^{\prime-1}$ is a subpath of the contour of $\Delta$,
and $|p'|\ge2\max\{|v_1\kl|,|v_2\kl|\}$,
at least one of the words $v_1^{\pm1}$, $v_2^{\pm1}$
is a subword of the label of $p'$, and therefore, it is
a subword of~$u_{i,j}^\sigma$.
Contradiction with condition~$(5)$.
Thus, $H$ is freely generated by $\phi(v_1\kl)$ and $\phi(v_2\kl)$.
The proof is complete.
\qed
\enddemo

\bigskip
The author would like to thank Alexander Ol'shanskii for suggesting the
problems and a possible approach to their solution.
\par

\enddocument